\documentclass[aap,seceqn,MSNbibl,dvips]{arximspdf}
\usepackage{mathbh}
\usepackage{accents}
\usepackage{graphicx}

\doi{10.1214/09-AAP676}
\volume{20}
\issue{6}
\pubyear{2010}
\firstpage{2040}
\lastpage{2085}

\makeatletter

\newcommand{\rme}{{e}}

\newtheorem{Theorem}{Theorem}[section]
\newtheorem{Proposition}[Theorem]{Proposition}
\newtheorem{Lemma}[Theorem]{Lemma}
\newproclaim{Remark}[Theorem]{Remark}

\makeatother

\begin{document}
\begin{frontmatter}

\title{Load optimization in a planar network}
\runtitle{Load optimization in a planar network}

\begin{aug}
\author[A]{\fnms{Charles} \snm{Bordenave}\corref{}\ead[label=e1]{charles.bordenave@math.univ-toulouse.fr}} and
\author[B]{\fnms{Giovanni Luca} \snm{Torrisi}\ead[label=e2]{torrisi@iac.rm.cnr.it}}
\runauthor{C. Bordenave and G. L. Torrisi}
\affiliation{CNRS, Universit\'e de Toulouse and CNR, Istituto per le
Applicazioni del~Calcolo~``Mauro~Picone''}
\address[A]{CNRS, Universit\'e de Toulouse\\
Institut de Math\'ematiques \\
118 route de Narbonne \\
31062 Toulouse\\
France \\
\printead{e1}}
\address[B]{CNR, Istituto per le Applicazioni \\
del Calcolo ``Mauro Picone''\\
c/o Department of Mathematics\\
University of Rome ``Tor Vergata'' \\
Via della Ricerca Scientifica 1\\
I-00133 Roma\\
Italia \\
\printead{e2}}
\end{aug}

\received{\smonth{2} \syear{2009}}
\revised{\smonth{12} \syear{2009}}

%
\begin{abstract}
We analyze the asymptotic properties of a Euclidean optimization
problem on the plane. Specifically, we consider a network with three
bins and $n$ objects spatially uniformly distributed, each object being
allocated to a bin at a cost depending on its position. Two allocations
are considered: the allocation minimizing the bin loads and the
allocation allocating each object to its less costly bin. We analyze
the asymptotic properties of these allocations as the number of objects
grows to infinity. Using the symmetries of the problem, we derive a law
of large numbers, a central limit theorem and a large deviation
principle for both loads with explicit expressions. In particular, we
prove that the two allocations satisfy the same law of large numbers,
but they do not have the same asymptotic fluctuations and rate
functions.
\end{abstract}

%
\begin{keyword}[class=AMS]
\kwd[Primary ]{60F05}
\kwd{60F10}
\kwd[; secondary ]{90B18}
\kwd{90C27}.
\end{keyword}
\begin{keyword}
\kwd{Euclidean optimization}
\kwd{law of large numbers}
\kwd{central limit theorem}
\kwd{large deviations}
\kwd{calculus of variations}
\kwd{wireless networks}.
\end{keyword}

\end{frontmatter}

\section{Introduction}\label{sec:int}

\begin{figure}

\includegraphics{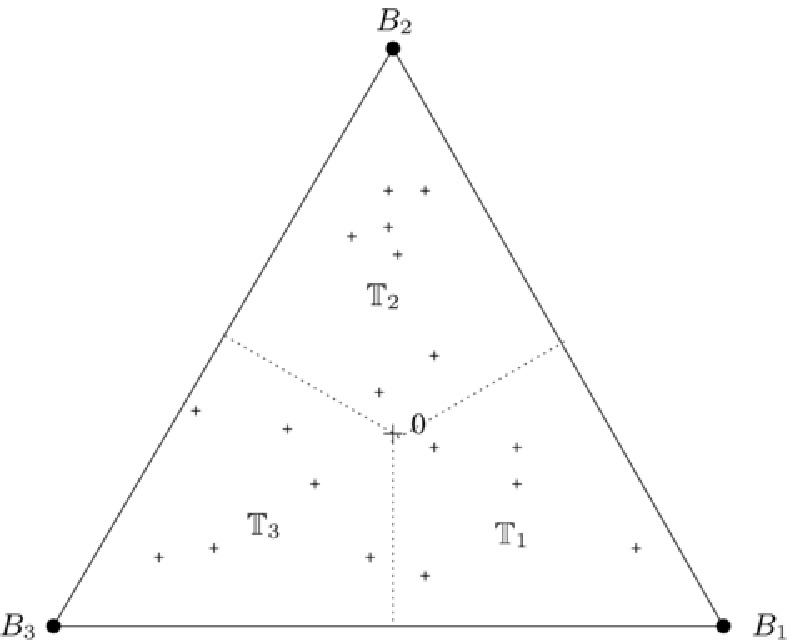}

\caption{The triangle $\mathbb{T}$, the three bins and the
$n$ objects.}\label{fig:T}
\end{figure}

In this paper we take an interest in a Euclidean optimization problem
on the
plane. For ease of notation, we shall identify the plane with the
set of complex numbers $\mathbb{C}$. Set $\lambda=2(3\sqrt3)^{-1/2}$,
$i=\sqrt{-1}$ (the complex unit), $j=\rme^{2i\pi/3}$ and
consider the triangle $\mathbb{T}\subset\mathbb{C}$ with\vspace*{1pt} vertices
$B_2=\lambda
i$, $B_1=j^2 B_2$ and $B_3=j B_2$. Note that $\mathbb{T}$ is an
equilateral triangle with side length $\lambda\sqrt{3}$ and unit
area. We label by $\{1,\ldots,n\}$ $n$ objects located in the
interior of $\mathbb{T}$ and denote by $X_k$, $k=1,\ldots,n$, the
location of the $k$th object; see Figure \ref{fig:T}. We assume that
$\{X_k\}_{k=1,\ldots,n}$ are independent random variables (r.v.'s)
with uniform distribution on $\mathbb{T}$. Suppose that there are three
bins located at each of the vertices of $\mathbb{T}$ and that each object
has to be allocated to a bin. The cost of an allocation is
described by a measurable function $c\dvtx\mathbb{T}\to[0,\infty)$
such that
$\|c\|_{\infty}:=\sup_{x\in\mathbb{T}}c(x)<\infty$. More precisely,
$c(x)=c_{1}(x)$ denotes the cost to allocate an object at
$x\in\mathbb{T}$ to the bin in~$B_1$; the cost to allocate an object at
$x\in\mathbb{T}$ to the bin in $B_2$ is $c_{2}(x)=c(j^{2}x)$; the cost
to allocate an object at $x\in\mathbb{T}$ to the bin in $B_3$ is
$c_{3}(x)=c(jx)$. Let
\[
\mathcal{A}_{n}=\bigl\{A=(a_{kl})_{1\leq k \leq n,1\leq l\leq3}\dvtx
a_{kl}\in\{0,1\}, a_{k1}+a_{k2}+a_{k3}=1\bigr\}
\]
be the set of allocation matrices: if $a_{kl}=1$ the
$k$th object is affiliated to the bin in $B_l$. We consider the
load relative to the allocation matrix $A=(a_{k l})_{1\leq k\leq
n,1\leq l\leq3}\in\mathcal{A}_n$:
\[
\rho_n(A)=\max_{1\leq l\leq3} \Biggl(\sum_{k=1}^{n}a_{kl}c_l
(X_k) \Biggr),
\]
and the minimal load
\[
\rho_n=\min_{A\in\mathcal{A}_n}\rho_n(A).
\]
Throughout this paper we refer to $\rho_n$ as the
optimal load. This simple instance of Euclidean optimization
problem has potential applications in operations research and
wireless communication networks. Consider three processors running
in parallel and sharing a pool of tasks $\{1,\ldots,n\}$ located,
respectively, at $\{X_1,\ldots,X_n\}\subseteq\mathbb{T}$. Suppose that
$c_l(x)$ is the time requested by the $l$th processor to process
a job located at $x\in\mathbb{T}$. Then $\rho_n$ is the minimal time
requested to process all jobs. For example, a natural choice
for the cost function is $c(x)=2|x-B_1|$, that is, the time of a
round-trip from $B_1$ to $x$ at unit speed. In a wireless
communication scenario, the bins are base stations and the objects
are users located at $\{X_1,\ldots,X_n\}\subseteq\mathbb{T}$. For the
base station located at $B_l$, the time needed to send one bit of
information to a user located at $x\in\mathbb{T}$ is $c_l(x)$. In this
context $\rho_n$ is the minimal time requested to send one bit of
information to each user and $1/\rho_n$ is the maximal throughput
that can be achieved. We have chosen a triangle $\mathbb{T}$ because it
is the fundamental domain of the hexagonal grid, which is a good model for
cellular wireless networks.

%
%

For $1\leq l\leq3$, we define the Voronoi cell associated to the
bin at $B_l$ by
\[
\mathbb{T}_{l}=\Bigl\{x\in\mathbb{T}\dvtx|x-B_l|={\min_{1\leq m\leq
3}}|x-B_m|\Bigr\}
\bigm\backslash D_l,
\]
where $D_1=\{ijt\dvtx t < 0\}$ and, for $l=2,3$, $D_l=\{ij^{l}t\dvtx t
\leq0\}$. Note that $\mathbb{T}_{1}\cup\mathbb{T}_{2}\cup\mathbb
{T}_{3}=\mathbb{T}$ and
$\mathbb{T}_{1}\cap\mathbb{T}_{2}=\mathbb{T}_{1}\cap\mathbb
{T}_{3}=\mathbb{T}_{2}\cap\mathbb{T}_{3}=
\varnothing$, that is, $\{\mathbb{T}_{1},\mathbb{T}_{2},\mathbb
{T}_{3}\}$ is a partition of
$\mathbb{T}$. Note also that $0\in\mathbb{T}_1$.

Throughout the paper, we denote by \mbox{$|\cdot|$} the Euclidean norm on
$\mathbb{C}$, by $\ell$ the Lebesgue measure on $\mathbb{C}$ and by
$x\cdot z$ the
usual scalar product on $\mathbb{C}$, that is, $x\cdot
z=\Re(x)\Re(z)+\Im(x)\Im(z)$. We suppose that the value of the
cost function is related to the distance of a point from a bin as
follows:
%
%
\begin{equation}\label{ass3}\hspace*{28pt}
\mbox{For all }x\in\mathbb{T}\mbox{ and }l=2,3\mbox{ if
}|x-B_1|<|x-B_l|\qquad\mbox{then } c_1(x)<c_l(x).
\end{equation}
For example, if $c(x)= f(|x-B_1|)$ and $f\dvtx[0,\infty)\to
[0,\infty)$ is increasing, then (\ref{ass3}) is satisfied.

In this paper, as $n$ goes to infinity, we study the properties of
an allocation which realizes the optimal load $\rho_n$, and, as a
benchmark, we compare it with the suboptimal load
$\overline{\rho}_n=\rho_n(\overline{A})$, where
$\overline{A}=(\overline{a}_{kl})_{1\leq k\leq n,1\leq l\leq3}$
is the random matrix obtained by affiliating each object to its
least costly bin
\[
\overline{a}_{kl}=\mathbh{1}(X_{k}\in\mathbb{T}_l).
\]
We shall prove that, using the strong symmetries of the
system, it is possible to perform a fine analysis of the
asymptotic optimal load. It turns out that a law of large number
can be deduced for the optimal and suboptimal load. More
precisely, setting
\[
\gamma=\int_{\mathbb{T}_1}c(x) \,dx,
\]
we have the following theorem.
\begin{Theorem}\label{th:lln} Assume
(\ref{ass3}). Then, almost surely (a.s.),
\[
\lim_{n\to\infty}\frac{\rho_n}{n}=\lim_{n\to\infty}\frac{
\overline{\rho}_n}{n}=\gamma.
\]
\end{Theorem}

As a consequence, at the first order, the optimal and
the suboptimal load perform similarly.

The next result shows that, at the second order, the two loads
differ significantly. We first introduce an extra symmetry
assumption on $c$, namely, its symmetry with respect to the
straight line determined by the points $0$ and $B_1$. If
$x=t\rme^{i\theta}\in\mathbb{T}$, $t
>0$, $\theta\in[0,2\pi]$, then its reflection with respect to the
straight line
determined by the points $0$ and $B_1$ is
$t\rme^{-i\theta-i{\pi}/{3}}\in\mathbb{T}$. Formally, we assume
%
%
\begin{eqnarray}\label{ass4}
&&c(t\rme^{i\theta})=
c(t\rme^{-i\theta-i{\pi}/{3}})\hspace*{29pt}\nonumber\\[-0pt]\\[-19pt]
\eqntext{\mbox{for all }\theta\in
[0,2\pi]\mbox{ and }t>0\mbox{ such that }t\rme^{i\theta} \in\mathbb{T}\mbox{ and}}\\
\eqntext{c\mbox{ is Lipschitz in a neighborhood of }D_{1}\cup
D_{3}.}
\end{eqnarray}
Setting
\[
\sigma^2=\int_{\mathbb{T}_1}c^{2}(x) \,dx
\]
and letting $\stackrel{d}{\rightarrow}$ denote the convergence in
distribution, we have the following theorem.
\begin{Theorem}\label{th:tcl}
Assume (\ref{ass3}) and (\ref{ass4}). Then, as $n$ goes to
infinity,
\[
n^{-1/2}(\rho_{n}-\gamma n)\stackrel{d}{\rightarrow}G,
\]
where $G$ is a Gaussian r.v. with zero mean and
variance $\sigma^{2}/3-\gamma^2$. Moroever, as $n$ goes to
infinity,
\[
n^{-1/2}(\overline\rho_{n}-\gamma
n)\stackrel{d}{\rightarrow}\max\{G_1,G_2,G_3\}-\tfrac
{1}{3}(G_{1}+G_{2}+G_{3}) + G
\]
and
\[
n^{-1/2}(\overline\rho_{n}-\rho_{n})\stackrel{d}{\rightarrow}\max
\{
G_1,G_2,G_3\}
-\tfrac{1}{3}(G_{1}+G_{2}+G_{3}),
\]
where $G_1$, $G_2$ and $G_3$ are independent Gaussian
r.v.'s with zero mean and variance~$\sigma^2$, independent of $G$. Finally
\[
\mathrm{E}[\rho_n]=n\gamma+o\bigl(\sqrt
n\bigr) \quad\mbox{and}\quad \mathrm{E}[\overline{\rho}_n]
=n\gamma+m\sqrt{n}+o\bigl(\sqrt n\bigr),
\]
where $m=\mathrm{E}[\max\{G_1,G_2,G_3\}]>0$ depends linearly on
$\sigma$.
\end{Theorem}

Theorem \ref{th:lln} states that $\overline\rho_n$ is
asymptotically optimal at scale $n$, but Theorem~\ref{th:tcl}
says that it is not asymptotically optimal at scale
$\sqrt{n}$. In the proof of Theorem \ref{th:tcl}, we shall
exhibit a suboptimal allocation which is asymptotically optimal at
scale $\sqrt{n}$ (see Proposition \ref{prop:subtcl}).

We shall also prove a large deviation principle (LDP) for both the
sequences $\{\rho_{n}/n\}_{n \geq1}$ and
$\{\overline{\rho}_{n}/n\}_{n \geq1}$. Recall that a family of
probability measures $\{\mu_n\}_{n\geq1}$ on a topological space
$(M,\mathcal{T}_M)$ satisfies a LDP with rate function $I$ if
$I\dvtx M\to[0,\infty]$ is a lower semi-continuous function such that
the following inequalities hold for every Borel set $B$
\[
-\inf_{y\in{\accentset{\circ}{B}}}
I(y)\leq\liminf_{n \rightarrow
\infty}\frac{1}{n}\log\mu_n (B)\leq\limsup
_{n\rightarrow\infty}\frac{1}{n}\log\mu_{n} (B)\leq-\inf_{y\in
\overline{B}}I(y),
\]
where $\accentset{\circ}{B}$ denotes the interior of $B$
and $\overline{B}$ denotes the closure of $B$. Similarly, we say
that a family of $M$-valued random variables $\{V_n\}_{n\geq1}$
satisfies an LDP if $\{\mu_n\}_{n\geq1}$ satisfies an LDP and
$\mu_n (\cdot)=P(V_n\in\cdot)$. We point out that the lower
semi-continuity of $I$ means that its level sets $\{y\in M\dvtx
I(y)\leq a\}$ are closed for all \mbox{$a\geq0$}; when the level sets
are compact the rate function $I(\cdot)$ is said to be good. For
more insight into large deviations theory, see, for instance, the
book by Dembo and Zeitouni \cite{dembo}.

We introduce an assumption on the level sets of the cost function
%
%
\begin{equation}\label{ass20}
\ell(c^{-1}(\{t\}))= 0\qquad\mbox{for all }t\geq0,
\end{equation}
an assumption on the regularity of $c$
%
%
\begin{equation}\label{ass2}
c\mbox{ is continuous on }\mathbb{T},
\end{equation}
and two further geometric conditions
%
%
\begin{eqnarray}\label{eq:concave}
&c(B_1)<c(x)<c(0)\qquad\mbox{for any }
x\in\mathbb{T}_1\setminus\{0,B_1\},&
\\
%
%
\label{eq:concave1}
&\displaystyle\frac{c_1(x)c_2(x)c_3(x)}{c_1(x)c_2(x)+c_1(x)c_3(
x)+c_2(x)c_3(x)}<\frac{c(0)}{3}<\int_{\mathbb{T}_2}c(z)
\,dz&\nonumber\\[-8pt]\\[-8pt]
\eqntext{\mbox{for any }x\in\mathbb{T}\setminus\{0\}.}
\end{eqnarray}
Assumption (\ref{eq:concave}) fixes the extrema of the cost
function on $\mathbb{T}_1$. The left-hand side inequality of
(\ref{eq:concave1}) imposes that $0$ is the most costly position
in terms of load [for a more precise statement, we postpone to (\ref
{eq:abcin})]. For $\theta\in\mathbb{R}$, define the functions
\[
\Lambda(\theta)=\log\biggl(3\int_{\mathbb{T}_1}\rme^{\theta
c(x)} \,dx \biggr) \quad\mbox{and}\quad
\overline{\Lambda}(\theta)=\log\biggl(\int_{\mathbb{T}_1}\rme
^{\theta
c(x)} \,dx+2/3 \biggr)
\]
and, for $y\in\mathbb{R}$, their Fenchel--Legendre transforms
\[
\Lambda^{*}(y)=\sup_{\theta\in\mathbb{R}}\bigl(\theta y-
\Lambda(\theta)\bigr) \quad\mbox{and}\quad \overline{\Lambda}{}^{*}(y)=
\sup_{\theta\in\mathbb{R}}\bigl(\theta y-\overline{\Lambda}(\theta)\bigr).
\]
The following LDPs hold:
\begin{Theorem}\label{th:LDP}
Assume (\ref{ass3}), (\ref{ass20}), (\ref{ass2}),
(\ref{eq:concave}) and (\ref{eq:concave1}). Then:

\begin{longlist}
\item $\{\rho_n/n\}_{n\geq1}$ satisfies an LDP on
$\mathbb{R}$ with good rate function
%
%
\begin{equation}
J(y)=\cases{
\Lambda^{*}(3y), &\quad if $y\in
\bigl(c(B_1)/3,c(0)/3\bigr)$,\cr
+\infty, &\quad otherwise.}
\end{equation}

\item $\{\overline{\rho}_n/n\}_{n\geq1}$ satisfies an
LDP on $\mathbb{R}$ with good rate function
%
%
\begin{equation}\label{eq:ratef}
\overline{J}(y)=\cases{
\Lambda^{*}(3y), &\quad if $y\in\bigl(c(B_1)/3,\gamma\bigr]$,\cr
\overline{\Lambda}{}^{*}(y), &\quad if $y\in(\gamma,c(0))$,\cr
+\infty, &\quad otherwise.}
\end{equation}
\end{longlist}
\end{Theorem}

The next proposition gives a more explicit expression
for the rate functions.
\begin{Proposition}\label{prop:Laplace}
Assume (\ref{ass3}), (\ref{eq:concave}) and $c$ continuous at
$0$ and $B_1$. Then $\Lambda^*$ and $\overline{\Lambda}{}^*$ are
continuous on
$(c(B_1),c(0))$ and

\[
\mbox{\textup{(i)} \hspace*{34pt}}\Lambda^{*}(y)=\cases{
y\theta_{y}-\Lambda(\theta_y), &\quad if $c(B_1)<y<c(0)$,\cr
+\infty, &\quad if $c(B_1)>y$ or $y>c(0)$,}\hspace*{34pt}
\]
where $\theta_y$ is the unique solution of
%
%
\begin{equation}\label{eq:tetay}
\frac{\int_{\mathbb{T}_1} c(x)\rme^{\theta
c(x)} \,dx}{\int_{\mathbb{T}_1}\rme^{\theta
c(x)} \,dx}=y;
\end{equation}
\[
\mbox{\textup{(ii)} \hspace*{36pt}}\overline{\Lambda}{}^{*}(y)=\cases{
y\eta_{y}-\overline{\Lambda}(\eta_y), &\quad if $c(B_1)<y<c(0)$,\cr
+\infty, &\quad if $c(B_1)>y$ or $y>c(0)$,}\hspace*{34pt}
\]
where $\eta_y$ is the unique solution of
%
%
\begin{equation}\label{eq:etay}
\frac{\int_{\mathbb{T}_1}c(x)\rme^{\theta
c(x)} \,dx}{\int_{\mathbb{T}_1}\rme^{\theta
c(x)} \,dx+2/3}=y.
\end{equation}
\item If $\gamma<y<c(0)/3$, then
$\overline{\Lambda}{}^{*}(y)<\Lambda^{*}(3y)$.
\end{Proposition}

\begin{figure}

\includegraphics{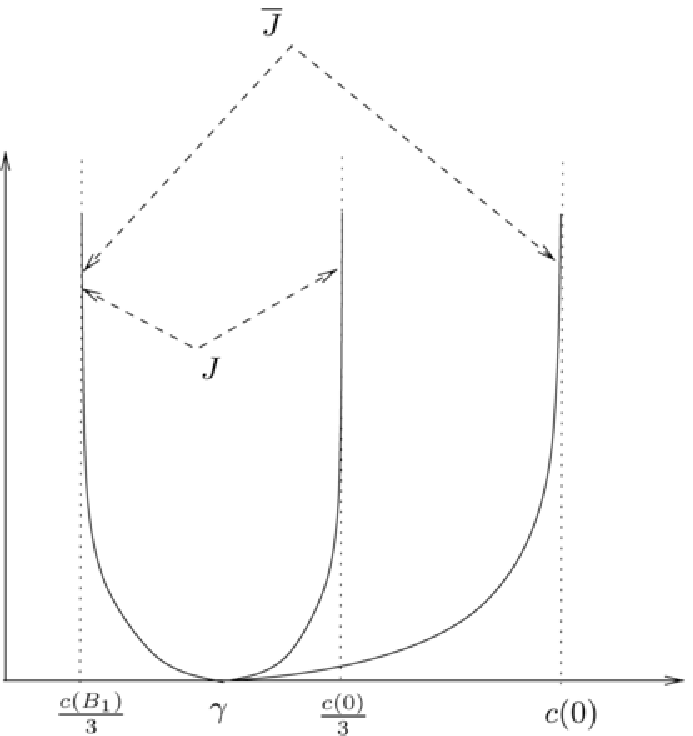}

\caption{The rate functions $J$ and $\overline{J}$.}\label{fig:rates}
\end{figure}

Note that $J(y)=\Lambda^*(3y)$ except possibly at
$y\in\{c(B_1),c(0)\}$; $\overline{J}(y)=\Lambda^*(3y)$ on
$(-\infty,\gamma]$ except possibly at $y=c(B_1)$, and
$\overline{J}(y)=\overline{\Lambda}{}^{*}(y)$ on $(\gamma,\infty)$ except
possibly at $y=c(0)$. These gaps are treated in Proposition
\ref{prop:extrema} with extra regularity assumptions on $c$. See
Figure \ref{fig:rates} for a schematic plot of the rate functions.
A simple consequence of Theorem \ref{th:LDP} and Proposition
\ref{prop:Laplace} is the following:
\[
\lim_{n\to\infty}\frac{\log P(\rho_{n}\geq n t)}{\log
P(\overline{\rho}_{n}\geq
nt)}=\frac{J(t)}{\overline{J}(t)} \quad\mbox{and}\quad \lim_{n\to\infty
}\frac
{P(\rho_{n}\geq
nt)}{P(\overline{\rho}_{n}\geq nt)}=0 \qquad \forall t\in\bigl(\gamma,c(0)/3\bigr).
\]
In words, it means that the probability of an
exceptionally large optimal load is significantly lower than the
probability of an exceptionally large suboptimal load; although,
on a logarithmic scale, the probability of an exceptionally small
optimal load does not differ significantly on the probability of
an exceptionally small suboptimal load. It is not in the scope of
this paper to discuss the trade-off between algorithmic complexity
and asymptotic performance. Moreover, we do not know if the
allocation that is asymptotically optimal at scale $\sqrt{n}$
used in the proof of Theorem \ref{th:tcl} (see Proposition
\ref{prop:subtcl}) has the same rate function than $\rho_n/n$.

%
%

Unlike it may appear, we shall not prove Theorem \ref{th:LDP} by
first computing the Laplace transform of $\rho_n$ and
$\overline{\rho}_n$ and then applying the G\"artner--Ellis theorem
(see, e.g., Theorem 2.3.6 in \cite{dembo}). We shall follow another
route. First, we combine Sanov's theorem (see, e.g., Theorem 6.2.10
in \cite{dembo}) and the contraction principle (see, e.g., Theorem
4.2.1 in \cite{dembo}) to prove that the sequences
$\{\rho_n/n\}_{n\geq1}$ and $\{\overline{\rho}_n/n\}_{n\geq1}$
obey a LDP, with rate functions given in variational form. Then,
we provide the explicit expression of the rate functions solving
the related variational problems. It is worthwhile to remark that,
using Theorem \ref{th:LDP} and Varadhan's lemma (see, e.g., Theorem
4.3.1 in \cite{dembo}) it is easily seen that
\[
\lim_{n\to\infty}\frac{1}{n}\log\mathrm{E}[\rme^{\theta\rho
_n}]=J^*(\theta) \quad\mbox
{and}\quad \lim_{n\to\infty}\frac{1}{n}\log\mathrm{E}[\rme
^{\theta
\overline
{\rho}_n}]=\overline{J}{}^*(\theta)\qquad \forall\theta\in\mathbb{R},
\]
where $J^*$ and $\overline{J}{}^*$ are the Fenchel--Legendre
transforms of $J$ and $\overline{J}$, respectively. A nice consequence of
Theorems \ref{th:lln} and \ref{th:tcl} is that, in terms of law of
the large numbers and central limit theorem, $\rho_n$ has the same
asymptotic behavior as
\[
\breve{\rho}_n=\frac{1}{3}\sum_{l=1}^{3}\sum_{k=1}^{n}\mathbh{1}\{
X_{k}\in\mathbb{T}
_l\}
c_l(X_k).
\]
Moreover, if the cost function satisfies extra regularity
assumptions (see Proposition \ref{prop:extrema}), by Theorem
\ref{th:LDP} and the G{\"a}rtner--Ellis theorem, we have that
$\rho_n$ and $\breve{\rho}_n$ have the same asymptotic behavior
even in terms of large deviations.

As can be seen from the proofs, if the left-hand side of
assumption (\ref{eq:concave1}) does not hold, then we have an
explicit rate function $J(y)$ only for $y<c(0)/3$. If the right-hand
side of assumption (\ref{eq:concave1}) also fails to hold,
then we have an explicit rate function $J(y)$ only for $y<y_0$ for
some $y_0>\gamma$. We also point out that the statements of
Theorems \ref{th:tcl} and \ref{th:LDP} concerning $\overline\rho_{n}$
do not require the use of (\ref{ass4}) and~(\ref{eq:concave}).

In wireless communication, the typical cost function is the
inverse of signal to noise plus interference ratio (see, e.g.,
Chapter IV in Tse and Viswanath \cite{tse}), which has the
following shape:
\[
c(x)=\frac{a+\min\{b,|x-B_2|^{-\alpha}\}+\min\{b,|x-
B_3|^{-\alpha}\}}{\min\{b,|x-B_1|^{-\alpha}\}},\qquad x\in\mathbb{T},
\]
where $\alpha\geq2$, $a>0$ and $b>(\lambda\sqrt{3}/2)^{-\alpha}$
[recall that $\lambda=2(3\sqrt{3})^{-1/2}$ and
$\lambda\sqrt{3}=|B_1-B_2|$]. We shall check in the
\hyperref[sec:app]{Appendix} that this cost function satisfies
(\ref{ass3}), (\ref{ass4}), (\ref{ass20}), (\ref{ass2}) and
(\ref{eq:concave}). Moreover, the first inequality in
(\ref{eq:concave1}) will be checked numerically and, for arbitrarily
fixed $\alpha>2$ and $a>0$, we shall determine values of the parameter
$b>(\lambda\sqrt{3}/2)^{-\alpha}$ such that the second inequality in
(\ref{eq:concave1}) holds.

The remainder of the paper is organized as follows. In Section
\ref{sec:lln} we analyze the sample path properties of the optimal
allocation and we prove Theorem \ref{th:lln}. In Section
\ref{sec:tcl} we show Theorem \ref{th:tcl}. Section \ref{sec:LDP}
is devoted to the proof of Theorem \ref{th:LDP} and Proposition
\ref{prop:Laplace}. In Section \ref{sec:extension}, we discuss
some generalizations of the model. We include also an \hyperref[sec:app]{Appendix}
where we prove some technical lemmas and provide an illustrative
example.

\section{Sample path properties}
\label{sec:lln}

\subsection{Structural properties of the optimal allocation}

Throughout this paper we denote by $\mathcal{M}_{b}(\mathbb{T})$ the
space of
Borel measures on $\mathbb{T}$ with total mass less than or equal to $1$
and by $\mathcal{M}_{1}(\mathbb{T})$ the space of probability
measures on
$\mathbb{T}$. These spaces are both equipped with the topology of weak
convergence (see, e.g., Billingsley \cite{billingsley}). For a
Borel function $h$ and a Borel measure $\mu$ on $\mathbb{T}$, we set
$\mu(h)=\int_{\mathbb{T}}h(x) \mu(dx)$. Consider the functional
from $\mathcal{M}_{b}(\mathbb{T})^3$ to $\mathbb{R}$ defined by
%
%
\begin{equation}\label{eq:littlephi}
\phi(\alpha_1,\alpha_2,\alpha_3)=\max
(\alpha_{1}(c_1),\alpha_2(c_2), \alpha_3(c_3)).
\end{equation}
Letting $\alpha_{|B}$ denote the restriction of a measure $\alpha$
to a Borel set $B$, we define the functionals $\Phi$ and $\Psi$
from $\mathcal{M}_{1}(\mathbb{T})$ to $\mathbb{R}$ by
\[
\Phi(\alpha)=\inf_{(\alpha_l)_{1\leq l\leq3}
\in\mathcal{M}_{b}(\mathbb{T})^3\dvtx\alpha_{1}+\alpha_{2}+\alpha_{3}
=\alpha}\phi(\alpha_1,\alpha_2,\alpha_3)
\]
and
\[
\Psi(\alpha)=\phi(\alpha_{|\mathbb{T}_1},\alpha_{|\mathbb
{T}_2},\alpha_{|\mathbb{T}_3}).
\]
Note that if $\delta_x$ denotes the Dirac measure with total mass
at $x\in\mathbb{T}$, then
%
%
\begin{equation} \label{eq:phirho}
\frac{\overline\rho_n}{n}=\Psi\Biggl(\frac{1}{n}\sum_{k=1}^{n}
\delta_{X_k} \Biggr).
\end{equation}
\begin{Lemma}\label{le:tvd}
Under assumption (\ref{ass2}) we have that $\phi$ is continuous
on $\mathcal{M}_{b}(\mathbb{T})^3$ and $\Psi$ and $\Phi$ are
continuous on
$\mathcal{M}_{1}(\mathbb{T})$ $($for the topology of the weak convergence$)$.
\end{Lemma}

The proof of Lemma \ref{le:tvd} is postponed to the \hyperref[sec:app]{Appendix}; the
continuity of $\phi$ and $\Psi$ is essentially trivial, but the
continuity of $\Phi$ requires more work. Define the set of
matrices
\[
\mathcal{B}_n=\{B=(b_{kl})_{1\leq k\leq n,1\leq l\leq3}\dvtx
b_{kl}\in[0,1], b_{k1}+b_{k2}+b_{k3}=1\}
\]
and
\[
\widetilde\rho_n=\min_{B\in\mathcal{B}_n}\rho_n(B).
\]
From the viewpoint of linear programming, this is the fractional
relaxation of the original optimization problem. Now, given a
matrix $B=(b_{kl})\in\mathcal{B}_n$, we define the associated measures
$(\alpha_1,\alpha_2,\alpha_3)\in\mathcal{M}_{b}(\mathbb{T})^3$ by
setting $\alpha_{l}=(1/n)\sum_{k=1}^{n}b_{kl}\delta_{X_k}$
($l=1,2,3$). Due to this correspondence, it is straightforward to
check that
%
%
\begin{equation}\label{eq:rhophi}
\frac{\widetilde\rho_n}{n}=\Phi\Biggl(\frac{1}{n}\sum
_{k=1}^{n}\delta
_{X_k} \Biggr).
\end{equation}

The next lemma is a collection of elementary statements, whose
proofs are given in the \hyperref[sec:app]{Appendix}.
\begin{Lemma}\label{le:tilderho}
Fix $n\geq1$ and let $B^{*}=(b^*_{kl})\in\mathcal{B}_n$ be an optimal
allocation matrix for $\widetilde{\rho}_n$. Then:
\begin{longlist}
\item For all $\alpha\in\mathcal{M}_1(\mathbb{T})$, there exists
$(\alpha_1,\alpha_2,\alpha_3)\in\mathcal{M}_b (\mathbb{T})^3$
such that $\alpha=
\alpha_{1}+\alpha_{2}+\alpha_{3}$ and
$\Phi(\alpha)=\phi(\alpha_1,\alpha_2,\alpha_3)$. Moreover,
whenever such equality holds, we have $\alpha_1(c_1)=\alpha_2(c_2)
=\alpha_3(c_3)$. In particular, the choice
$\alpha_l=(1/n)\sum_{k=1}^{n}b^*_{kl}\delta_{X_k}$ $(l=1,2,3)$
yields
\[
\sum_{k=1}^{n}b^{*}_{k1}c_{1}(X_k)=\sum_{k=1}^{n}b^{*}_{k2}c_{2}
(X_k)=\sum_{k=1}^{n}b^{*}_{k3}c_{3}(X_k).
\]
\item If assumption (\ref{ass20}) holds, then
\[
\rho_{n}-3\|c\|_\infty\leq
\widetilde{\rho}_n\leq\rho_n \qquad\mbox{a.s.}
\]
\item If assumption (\ref{ass20}) holds then the
sequences $\{\widetilde\rho_n/n\}$ and $\{\rho_n/n\}$ are
exponentially equivalent.
\end{longlist}
\end{Lemma}

For the definition of exponential equivalence, see page 130 in
\cite{dembo}.

\subsection{\texorpdfstring{Proof of Theorem \protect\ref{th:lln}}{Proof of Theorem 1.1}}

The law of large numbers yields, for all $l=1,2,3$,
\[
\lim_{n\to\infty}\frac{1}{n}\sum_{k=1}^{n}c_{l}(X_k)\mathbh{1}\{
X_k\in
\mathbb{T}_l\}=\int_{\mathbb{T}_l}c_{l}(x) \,dx=\gamma\qquad\mbox{a.s.}
\]
Therefore from the identity
\[
\frac{\overline{\rho}_n}{n}=\max_{1\leq l\leq3}\frac{1}{n}
\sum_{k=1}^{n}c_{l}(X_k)\mathbh{1}\{X_k\in\mathbb{T}_l\},
\]
we get $\lim_{n\to\infty}\overline{\rho}_n/n=\gamma$ a.s. We also
have to prove that $\lim_{n\to\infty} \rho_n/n=\gamma$ a.s. Let
$A=(a_{kl})\in\mathcal{A}_n$ be an allocation matrix. By assumption
(\ref{ass3}), if $x\in\mathbb{T}_l$ then $c_l(x)=\min_{1\leq m\leq3}
c_m (x)$. Therefore
%
%
\begin{eqnarray}\label{eq:3lln}
3\rho_{n}(A) & \geq& \sum_{l=1}^{3}\sum_{k=1}^{n}a_{kl }c_{l}(X_k)
\nonumber\\
& \geq& \sum_{l=1}^{3}\sum_{X_{k}\in\mathbb{T}_l}c_{l}(X_k)
\\
& \geq& 3\min_{1\leq l\leq3}
\Biggl(\sum_{k=1}^{n}c_{l}(X_k)\mathbh{1}\{X_k\in\mathbb{T}_l\}\Biggr).\nonumber
\end{eqnarray}
So taking the minimum over all the allocation matrices
we deduce
\[
\min_{1\leq l\leq3} \Biggl(\sum_{k=1}^{n}c_{l}(X_k)\mathbh{1}\{X_k\in
\mathbb{T}_l\} \Biggr)\leq\rho_{n}\leq\overline{\rho}_n.
\]
Thus by applying the law of large numbers, we have a.s.
\[
\gamma\leq\liminf_{n\to\infty}\frac{\rho_n}{n}\leq\limsup
_{n\to\infty
}\frac{\rho_n}{n}\leq
\gamma.
\]

\begin{Remark}\label{remark1}
Assume that conditions (\ref{ass3}), (\ref{ass20})
and (\ref{ass2}) hold. By Theorem~\ref{th:lln} we have
$\lim_{n\to\infty}\overline{\rho}_n/n=\gamma$ a.s. So by Lemma
\ref{le:tvd}, equation (\ref{eq:phirho}) and the a.s. weak
convergence of $(1/n)\sum_{k=1}^{n}\delta_{X_k}$ to $\ell$ we get
$\Psi(\ell)=\gamma$. Similarly, using (\ref{eq:rhophi})
in place of (\ref{eq:phirho}), we deduce that
$\lim_{n\to\infty}\widetilde{\rho}_n/n=\Phi(\ell)$ a.s. By
Lem\-ma~\ref{le:tilderho}\textup{(ii)}, $| \widetilde{\rho}_n/n -\rho_n /n |
\leq
3 \| c\|_\infty/ n$, so we obtain
$\lim_{n\to\infty}\rho_n/n=\Phi(\ell)$ a.s., and by Theorem
\ref{th:lln} we have $\Phi(\ell)=\gamma$.
\end{Remark}

\section{\texorpdfstring{Proof of Theorem \protect\ref{th:tcl}}{Proof of Theorem 1.2}}
\label{sec:tcl}

Consider the random signed measure
\[
W_{n}=\sqrt{n} (\mu_n-\ell) \qquad\mbox{where } \mu_n=
\frac1 n \sum_{k=1}^{n}\delta_{X_k}.
\]
The standard Brownian bridge $W$ on $\mathbb{T}$ is a random signed
measure specified by the centered Gaussian process $\{W(f)\}$
(indexed on the set of square integrable functions on $\mathbb{T}$, with
respect to $\ell$), with covariance given by
\[
\mathrm{E}[W(f)W(g)]=\ell(fg)-\ell(f)\ell(g)
\]
(see, e.g., Dudley \cite{dudley}). By construction,
\[
\overline{\rho}_{n}=n\max_{1\leq l\leq3} \biggl(\int_{\mathbb{T}_{l}}
c_l(x)\mu_n(dx) \biggr)
\]
or equivalently
%
%
\begin{equation} \label{eq:tcla}
\frac{\overline{\rho}_{n}-n\gamma}{\sqrt n}=\max_{1\leq l\leq3}
\biggl(\int_{\mathbb{T}_l}c_l(x)W_n(dx) \biggr).
\end{equation}
Let $f$ be a square integrable function on $\mathbb{T}$. Then,
as $n\to\infty$,
\[
W_n(f)=\frac{\sum_{k=1}^{n}f(X_k)-n\ell(f)}{\sqrt{n}}
\stackrel{d}{\rightarrow}W(f).
\]
Indeed, by the central limit theorem $W_n(f)$ converges
in distribution to a Gaussian r.v. with zero mean and variance
equal to $\ell(f^2)-\ell^2(f)$, which is exactly the law of
$W(f)$. Using the L{\'e}vy continuity theorem and the inversion
theorem (see, e.g., Theorems 7.5 and 7.6 in \cite{billingsley}), we
have, for all square integrable functions $f_1$, $f_2$ and $f_3$,
\[
(W_n(f_1),W_n(f_2),W_n(f_3))\stackrel{d}{\rightarrow}
(W(f_1),W(f_2),W(f_3)).
\]
For $(x_1,x_2,x_3)\in\mathbb{R}^3$, the function $(x_1,x_2,x_3)
\mapsto\max(x_1, x_2,x_3)$ is continuous. Therefore, by the
continuous mapping theorem (see, e.g., Theorem 5.1 in
\cite{billingsley}) and (\ref{eq:tcla}) we have, as $n$ goes to
infinity,
%
%
\begin{equation}\label{eq:convdibbar}
\frac{\overline{\rho}_{n}-n\gamma}{\sqrt
n}\stackrel{d}{\rightarrow}\max_{1\leq l\leq
3} \biggl(\int_{\mathbb{T}_l}c_l(x)W(dx) \biggr).
\end{equation}
We shall show later on that the r.v. in the right-hand
side of (\ref{eq:convdibbar}) has the claimed distribution. Now
we consider the optimal load $\rho_n$. By the second inequality in
(\ref{eq:3lln}) we have
\[
3\rho_{n}\geq
n\sum_{l=1}^{3}\int_{\mathbb{T}_l}c_{l}(x)\mu_{n}(dx)
\]
and therefore
%
%
\begin{equation}\label{eq:modlb}
3\frac{\rho_{n}-n\gamma}{\sqrt{n}}\geq\sum_{l=1}^{3}
\int_{\mathbb{T}_l}c_{l}(x)W_{n}(dx).
\end{equation}
The following proposition is the heart of the proof. It will be
shown later on.
\begin{Proposition}\label{prop:subtcl}
Under the assumptions of Theorem \ref{th:tcl}, there exist
absolute constants $L_0$ and $L_1$, not depending on $n$, such
that the following holds. For any $1/4<\alpha<1/2$, with
probability at least $1-L_{1}\exp(-L_{0}n^{1-2\alpha})$, there
exists an allocation matrix $\hat{A}=(\hat{a}_{kl})_{1\leq k\leq
n,1\leq l\leq3}\in\mathcal{A}_n$ with associated load
$\hat{\rho}_{n}= \rho_n(\hat{A})$ such that
\[
\Biggl|3\frac{\hat{\rho}_{n}-n\gamma}{\sqrt{n}}-\sum_{l=1}^{3}
\int_{\mathbb{T}_l}c_l(x)W_n(dx) \Biggr|\leq n^{1/2-2\alpha}.
\]
\end{Proposition}

Using this result, $\hat{\rho}_{n}\geq\rho_n$ and
(\ref{eq:modlb}), we have that with probability at least
$1-L_{1}\exp(-L_{0}n^{1-2\alpha})$
%
%
\begin{equation}\label{eq:tcl10}
\Biggl|3\frac{\rho_{n}-n\gamma}{\sqrt{n}}-\sum_{l=1}^{3}
\int_{\mathbb{T}_l}c_l(x)W_n(dx) \Biggr|\leq n^{1/2-2\alpha}.
\end{equation}
Therefore, as $n$ goes to infinity,
\[
\frac{\rho_{n}-n\gamma}{\sqrt{n}}-\frac{1}{3}\sum_{l=1}^{3}
\int_{\mathbb{T}_l}c_l(x)W_n(dx)\stackrel{d}{\rightarrow}0.
\]
The continuous mapping theorem yields
\[
\sum_{l=1}^{3}
\int_{\mathbb{T}_l}c_l(x)W_n(dx)\stackrel{d}{\rightarrow}\sum_{l=1}^{3}
\int_{\mathbb{T}_l}c_l(x)W(dx).
\]
So combining these latter two limits we get, as $n$ goes to
infinity,
\[
\frac{\rho_{n}-n\gamma}{\sqrt{n}}\stackrel{d}{\rightarrow}\frac{1}{3}
\sum_{l=1}^{3}\int_{\mathbb{T}_{l}}c_l(x)W(dx),
\]
that is, $n^{-1/2}(\rho_{n}-n\gamma)$ converges weakly to a
centered Gaussian random variable with variance
$\sigma^{2}/3-\gamma^2$. We have considered so far, the normalized
sequences $\rho_n$ and $\overline{\rho}_n$ separately. However, we
can carry the same analysis on the normalized difference
$\overline{\rho}_n-\rho_n$. More precisely, by
(\ref{eq:tcla}) we have a.s.
\begin{eqnarray*}
&& \Biggl|\frac{\overline{\rho}_{n}-\rho_n}{\sqrt{n}}- \Biggl[\max
_{1\leq
l\leq
3} \biggl(\int_{\mathbb{T}_l}c_l(x)W_n(dx) \biggr)-\frac{1}{3}
\sum_{l=1}^{3}\int_{\mathbb{T}_l}c_l(x)W_n(dx) \Biggr]
\Biggr|
\\
&&\qquad\leq\biggl|\frac{\overline{\rho}_{n}-n\gamma}{\sqrt n}-\max
_{1\leq
l\leq3} \biggl(\int_{\mathbb{T}_l}c_l(x)W_n(dx) \biggr) \biggr|\\
&&\qquad\quad{} +
\Biggl|\frac{\rho_{n}-n\gamma}{\sqrt{n}}-\frac{1}{3}\sum_{l=1}^{3}
\int_{\mathbb{T}_l}c_l(x)W_n(dx) \Biggr|\\
&&\qquad= \Biggl|\frac{\rho_{n}-n\gamma}{\sqrt{n}}-\frac{1}{3}\sum_{l=1}^{3}
\int_{\mathbb{T}_l}c_l(x)W_n(dx) \Biggr|.
\end{eqnarray*}
Thus, by (\ref{eq:tcl10}), we obtain, with probability
at least $1-L_{1} \exp(-L_{0}n^{1-2\alpha})$,
\[
\Biggl|\frac{\overline{\rho}_{n}-\rho_n}{\sqrt{n}}- \Biggl[\max
_{1\leq
l\leq
3} \biggl(\int_{\mathbb{T}_l}c_l(x)W_n(dx) \biggr)-\frac{1}{3}
\sum_{l=1}^{3}\int_{\mathbb{T}_l}c_l(x)W_n(dx) \Biggr]
\Biggr|\leq
\frac{1}{3}n^{1/2-2\alpha}.
\]
Therefore, as $n\to\infty$,
\[
\frac{\overline{\rho}_{n}-\rho_n}{\sqrt{n}}- \Biggl[\max_{1\leq
l\leq
3} \biggl(\int_{\mathbb{T}_l}c_l(x)W_n(dx) \biggr)-\frac{1}{3}
\sum_{l=1}^{3}\int_{\mathbb{T}_l}c_l(x)W_n(dx) \Biggr]\stackrel
{d}{\rightarrow}0.
\]
The continuous mapping theorem yields
\begin{eqnarray*}
&&\max_{1\leq l\leq
3} \biggl(\int_{\mathbb{T}_l}c_l(x)W_n(dx) \biggr)-\frac{1}{3}
\sum_{l=1}^{3}\int_{\mathbb{T}_l}c_l(x)W_n(dx)\\
&&\qquad \stackrel{d}{\rightarrow}\max_{1\leq l\leq
3} \biggl(\int_{\mathbb{T}_l}c_l(x)W(dx) \biggr)-\frac{1}{3}
\sum_{l=1}^{3}\int_{\mathbb{T}_l}c_l(x)W(dx)
\end{eqnarray*}
and therefore, as $n\to\infty$,
\[
\frac{\overline{\rho}_{n}-\rho_n}{\sqrt
n}\stackrel{d}{\rightarrow}\max_{1\leq l \leq
3} \biggl(\int_{\mathbb{T}_l}c_l(x)W(dx) \biggr)-\frac{1}{3}
\sum_{l=1}^{3}\int_{\mathbb{T}_l}c_l(x)W(dx).
\]
For $l\in\{1,2,3\}$, set
\[
N_l=\int_{\mathbb{T}_l}c_l(x)W(dx)-\frac{1}{3}\sum_{l=1}^{3}\int_{
\mathbb{T}_l}c_l(x)W(dx).
\]
By definition $\{W(f)\}$ is a centered Gaussian process
indexed on the set of square integrable functions; therefore
$N=(N_1,N_2,N_3)$ follows a multivariate Gaussian distribution
with mean $0$. A simple computation shows that the covariance
matrix of $N$ is
\[
\frac{\sigma^2}{3} \pmatrix{
2 & -1 & -1 \cr
-1 & 2 & -1 \cr
-1 & -1 & 2}.
\]
It implies that $N$ has the same distribution as
\[
\bigl(G_1-(G_{1}+G_{2}+G_{3})/3,G_2-(G_{1}+G_{2}+G_{3})/3,G_{3}-(G_{1}
+G_{2}+G_{3})/3\bigr),
\]
where $G_1$, $G_2$ and $G_3$ are independent Gaussian r.v.'s with
mean $0$ and variance~$\sigma^2$. Moreover $N$ is independent of
$\frac
{1}{3}\sum_{l=1}^{3}\int_{
\mathbb{T}_l}c_l(x)W(dx)$, and we deduce the claimed expression for
(\ref{eq:convdibbar}).

It remains to compute the asymptotic behavior of the expectation
of the loads. A direct computation gives, for any $l=1,2,3$,
\[
\mathrm{E}\biggl[ \biggl(\int_{\mathbb{T}_l}c_l(x)W_n(dx) \biggr)^2 \biggr]
=\frac{\sigma^2}{3}-\frac{\gamma^2}{9n}\leq\frac{\sigma^2}{3}.
\]
Thus the sequences $\{\int_{\mathbb{T}_l}c_l(x)W_n(dx)\}$
$(l=1,2,3)$ are uniformly integrable. This implies that the
sequence $ \{\max_{1\leq l\leq
3} (\int_{\mathbb{T}_l}c_l(x)W_n(dx) ) \}$ is
uniformly integrable and so using (\ref{eq:tcla}) we
have
\begin{eqnarray*}
\lim_{n\to\infty}\mathrm{E}[\overline{\rho}_{n}-n\gamma]/
\sqrt{n}&=&\lim_{n\to\infty}\mathrm{E}\biggl[\max_{1\leq l\leq3}
\biggl(\int_{\mathbb{T}_l}c_l(x)W_n(dx) \biggr) \biggr]\\
&=&\mathrm{E}\biggl[\max_{1\leq l\leq3}
\biggl(\int_{\mathbb{T}_l}c_l(x)W(dx) \biggr) \biggr]\\
&=&m=\mathrm{E}[\max\{G_1,G_2,G_3\} ].
\end{eqnarray*}
Now we give the asymptotic behavior of
$\mathrm{E}[\rho_n]$. Note that by (\ref{eq:tcl10}) we have
\begin{eqnarray*}
&&\mathrm{E} \Biggl[ \Biggl|3\frac{\rho_{n}-n\gamma}{\sqrt{n}}-\sum_{l=1}^{3}
\int_{\mathbb{T}_l}c_l(x)W_n(dx) \Biggr| \Biggr]\\
&&\qquad\leq
n^{1/2-2\alpha}+\mathrm{E} \Biggl[ \Biggl|3\frac{\rho_{n}-n\gamma
}{\sqrt
{n}}-\sum_{l=1}^{3}
\int_{\mathbb{T}_l}c_l(x)W_n(dx) \Biggr|\mathbh{1}\{|\cdots|>
n^{1/2-2\alpha}\} \Biggr]\\
&&\qquad\leq n^{1/2-2\alpha}+10\|c\|_{\infty}L_{1}\sqrt{n}
\exp(-L_{0}n^{1-2\alpha})\\
&&\qquad=n^{1/2-2
\alpha}+\widetilde{L}_{1}\sqrt{n}\exp(-L_{0}n^{1-2\alpha}),
\end{eqnarray*}
where the latter inequality follows since $\gamma\leq
\|c\|_{\infty}$, $\rho_{n}\leq\|c\|_{\infty}n$ and
$|{\int_{\mathbb{T}_l}}c_l(x)\times\break W_n(dx)|\leq
2\|c\|_{\infty}\sqrt{n}$. Therefore,
since $\mathrm{E}[\int_{\mathbb{T}_{l}}c_l(x)W_n(dx) ]=0$
and $1/4
< \alpha<1/2$, our computation leads to
\[
\lim_{n\to\infty}\mathrm{E}[\rho_{n}-n\gamma]/\sqrt{n}=0.
\]
\begin{pf*}{Proof of Proposition \ref{prop:subtcl}}
We start describing the allocation matrix $\hat{A}$. For
$l,m\in\{1,2,3\}$ and $t\in[- \lambda\sqrt{3}/2,\lambda\sqrt{3}/2]$,
denote by $B_{lm}(t)$ the point on the segment $\overline{B_{l}B_{m}}$
at distance $t+\lambda\sqrt{3}/2$ from $B_l$. We extend the definition
of $B_{lm}(t) $ for all $t\in[-\lambda\sqrt{3},\lambda\sqrt{3}]$ by
following the edges of $\mathbb{T}$. More precisely, we set
\[
B_{12}(t)=\cases{
B_{31}\bigl(\lambda\sqrt{3}+t\bigr), &\quad if $t \in
\bigl[-\lambda\sqrt{3},-\lambda\sqrt{3}/2\bigr]$,\cr
B_{23}\bigl(\lambda\sqrt{3}-t\bigr), &\quad if $t \in
\bigl[\lambda\sqrt{3}/2,\lambda\sqrt{3}\bigr]$.}
\]
For $l,m\in\{1,2,3\}$, $B_{lm}(t)$ is defined similarly by a
circular permutation of the indices. For $\mathbf{t}=(t^1,t^2,t^3)\in[-
\lambda\sqrt{3},\lambda\sqrt{3}]^3$, let
\[
C^{1}(\mathbf{t})=\{0\}\cup\bigl(\{z\in\mathbb{C}\dvtx
z\cdot(B_{12}(t^1)\rme^{-i\pi/2})\geq0\}\cap\{z\in\mathbb
{C}\dvtx
z\cdot(B_{31}(t^3)\rme^{i\pi/2})>0\} \bigr)
\]
be the (possibly empty) cone delimited by the straight
line determined by the points~$0$, $B_{12}(t^1)$ and $B_{31}(t^3)$. We define
$\Gamma^{1}(\mathbf{t})=C^{1}(\mathbf{t})\cap\mathbb{T}$.
Similarly, let
$\Gamma^{2}(\mathbf{t})=C^{2}(\mathbf{t})\cap\mathbb{T}$ and
$\Gamma
^{3}(\mathbf{t})=C^{3}(\mathbf{t})\cap
\mathbb{T}$
with
\begin{eqnarray*}
C^{2}(\mathbf{t})&=&\{z\in\mathbb{C}\dvtx z\cdot(B_{12}(t ^1)\rme
^{i\pi
/2})>0\}
\cap\{z\in\mathbb{C}\dvtx z\cdot(B_{23}(t^2)\rme^{-i\pi/2})\geq
0\},\\
C^{3}(\mathbf{t})&=&\{z\in\mathbb{C}\dvtx z\cdot(B_{23}(t^2)\rme
^{i\pi
/2})>0\}
\cap\{z\in\mathbb{C}\dvtx z\cdot(B_{31}(t^3)\rme^{-i\pi/2})\geq
0\}.
\end{eqnarray*}
By construction, the sets $\Gamma^{1}(\mathbf{t})$, $\Gamma
^{2}(\mathbf{t})$ and
$\Gamma^{3}(\mathbf{t})$ are disjoint and their union is $\mathbb
{T}$. For
$l\in\{1,2,3\}$, set
\[
\rho_{n}^{l}(\mathbf{t})=\sum_{k=1}^{n}c_{l}(X_k)\mathbh{1}\{X_k\in
\Gamma^{l}(\mathbf{t})\}
\]
and consider the following recursion. At step $0$: for
$\mathbf{t}_0=(0,0,0)$, define
\[
m_{0}=\mathop{\arg\min}_{1\leq l\leq3}\rho_{n}^{l}(\mathbf{t}_0)
\]
(breaking ties with the lexicographic order) and
\[
M_{0}=\mathop{\arg\max}_{1\leq l\leq3}\rho_{n}^{l}(\mathbf{t}_0)
\]
(again breaking ties with the lexicographic order). If
$\rho_{n}^{M_0}(\mathbf{t}_0)-\rho_{n}^{m_0}(\mathbf{t}_0)\leq2\|
c\|_\infty$,
the recursion stops. Otherwise,
$\rho_{n}^{M_0}(\mathbf{t}_0)-\rho_{n}^{m_0}(\mathbf{t}_0)>2\|c\|
_\infty$ and
there is at least one point $X_i$ $(i=1,\ldots,n)$ in
$\Gamma^{M_0}(\mathbf{t}_0)$. Note also that, a.s., for all $\theta
\in[0,2
\pi]$, there is at most one point of $\{X_1,\ldots,X_n\}$ on the
straight line $(x\rme^{i\theta}, x >0)$. As a consequence
there exists a random variable $0\leq t_{1}\leq\lambda\sqrt{3}$
such that, a.s., there is exactly one point $X_i$ $(i=1,\ldots,n)$
in the triangle with vertices $\{0,B_{m_{0}
M_{0}}(t_1),B_{m_{0}M_{0}}(0)\}$ for $0\leq t_{1}\leq\lambda\sqrt
{3}/2$, or in the polygon with vertices $\{0,B_{m_{0}
M_{0}}(t_1),B_{M_0}, B_{m_{0}M_{0}}(0)\}$ for $\lambda\sqrt{3}/2<
t_{1}\leq\lambda\sqrt{3}$. We then set
$\mathbf{t}_{1}=(t_1^1,t_1^2,t_1^3):=(t_1,0,0)$ if $m_{0}=1$, $M_{0}=2$;
$\mathbf{t}_{1}=(- t_1,0,0)$ if $m_{0}=2$, $M_{0}=1$; $\mathbf
{t}_{1}=(0,t_1,0)$
if $m_{0} =2$, $M_{0}=3$; $\mathbf{t}_{1}=(0,-t_1,0)$ if $m_{0}=3$,
$M_{0}=2$; $\mathbf{t}_{1}=(0,0,-t_1)$ if $m_{0}=1$, $M_{0}=3$;
$\mathbf{t}_{1}=(0,0,t_1)$ if $m_{0}=3$, $M_{0}=1$. The sets $(\Gamma
^1(\mathbf{t}
_1),\Gamma^2(\mathbf{t}_1),\Gamma^3(\mathbf{t}_1))$ are thus
designed to allocate one
extra point to bin $m_0$ and one less to $M_0$. By construction, we
have
\[
\rho_{n}^{m_0}(\mathbf{t}_{1})<\rho_{n}^{M_0}(\mathbf{t}_{1})
,\qquad
\max_{1\leq l\leq3}\rho_{n}^{l}(\mathbf{t}_{1})<\max_{1\leq l\leq3}
\rho_{n}^{l}(\mathbf{t}_{0})
\]
and
\[
\min_{1\leq l\leq3}
\rho_{n}^{l}(\mathbf{t}_{1})>\min_{1\leq l\leq3}\rho
_{n}^{l}(\mathbf{t}_{0}).
\]
At step 1: define
\[
m_{1}=\mathop{\arg
\min}_{1\leq l\leq3}\rho_{n}^{l}(\mathbf{t}_1)
\]
(breaking ties with the lexicographic order) and
\[
M_1=\mathop{\arg\max}_{1\leq l\leq3}
\rho_{n}^{l}(\mathbf{t}_1)
\]
(again breaking ties with the lexicographic
order). Similarly to step $0$, if
$\rho_{n}^{M_1}(\mathbf{t}_{1})-\rho_{n}^{m_1}(\mathbf{t}_{1})>2\|
c\|_{\infty}$,
then there is at least one point of $\{X_1,\ldots,X_n\}$ in
$\Gamma^{M_1}(\mathbf{t}_1)$ and we build the random vector
$\mathbf{t}_2=(t_2^1,t_2^2,t_2^3)$ in order to allocate one extra
point to bin
$m_1$ and one less to $M_1$. The recursion stops at the first step
$k\geq0$ such that
\[
\rho_{n}^{M_{k}}(\mathbf{t}_{k})-\rho_{n}^{m_{k}}(\mathbf
{t}_{k})\leq
2\|c\|_\infty
\]
(where $m_{k}$, $M_{k}$ and $\mathbf{t}_{k}$ are defined similarly to
$m_0,m_1,\ldots, M_0,M_1,\ldots$ and $\mathbf{t}_1,\mathbf
{t}_2,\ldots$). As
we shall check soon, the recursion stops after at most $n$
steps. When the recursion stops, say at step $k_n\leq n$, we set
$\Gamma^{l}_{n}=\Gamma^{l}(\mathbf{t}_{k_n})$ and $\mathbf
{t}_{n}=\mathbf{t}_{k_n}$. The
allocation matrix $\hat{A}$ is defined by allocating $X_k$ to the
bin in $B_l$ if $X_k\in\Gamma^{l}_{n}$, that is,
\[
\hat{A}=(\hat{a}_{kl})_{1\leq k\leq n, 1\leq l\leq3}\qquad
\mbox{where }\hat{a}_{kl}=\mathbh{1}\{X_k\in\Gamma_n^l\}.
\]
By construction, we have for all $l,m \in\{1,2,3\}$,
%
%
\begin{equation}\label{eq:hatAload}
|\rho^{l}_{n}(\mathbf{t}_{n})-\rho^{m}_{n}(\mathbf{t}_{n})|\leq2\|
c\|_{\infty}.
\end{equation}
We now analyze the recursion more closely. Assume that
at step $0$ we have $m_0=3$ and $M_0=1$,
that is, $\rho^{1}_{n}(\mathbf{t}_0)\geq\rho^{2}_{n}(\mathbf
{t}_0)\geq\rho^{3}_{n}
(\mathbf{t}_0)$. Then, for all $k\leq k_n$,
%
%
\begin{equation}\label{eq:recineq1}
\rho^1_{n}(\mathbf{t}_k)\geq\rho^2_{n}(\mathbf{t}_k)-\|c\|_\infty
\quad\mbox{and}\quad \rho^3_{n}(\mathbf{t}_k)\leq\rho^2_{n}(\mathbf{t}_k)+\|c\|
_\infty.
\end{equation}
Indeed, if for all $k<k_n$, $m_{k}=3$ and $M_{k}=1$, there is
nothing to prove since $|\rho^l_{n}(\mathbf{t}_{k+1})-\rho^l_{n}
(\mathbf{t}_{k})|\leq\|c\|_{\infty}$. Assume that there exists $k<k_n$
such that $m_{k}\neq3$ or $M_{k}\neq1$. We define
\[
k_{0}=\min\{k\geq1\dvtx m_k \neq3 \mbox{ or } M_k \neq1 \}.
\]
For concreteness, assume, for example, that $M_{k_0} \neq
1$. By construction, $k_{0}-1<k_n$ so that
$\rho_{n}^{1}(\mathbf{t}_{k_{0}-1})
>\rho_{n}^{3}(\mathbf{t}_{k_{0}-1})+2\|c\|_\infty$. Since $\rho_n^1(
\mathbf{t}_{k_{0}-1})\geq\rho_{n}^{2}(\mathbf{t}_{k_{0}-1})\geq
\rho_{n}^{3}(
\mathbf{t}_{k_{0}-1})$, we deduce that $M_{k_0}=2$ and
$m_{k_0}=3$. Recall that, for $k<k_n$, $\rho_{n}^{M_{k}}(\mathbf{t}_{k})-
\|c\|_{\infty}\leq\rho_{n}^{M_{k}}(\mathbf{t}_{k+1})<
\rho_{n}^{M_{k}}(\mathbf{t}_{k})$. Thus, for $k=k_{0}-1$, from $ \rho^1
_{n}(\mathbf{t}_{k_0})\leq\rho^{2}_{n}(\mathbf{t}_{k_0})=\rho
^{2}_{n}(\mathbf{t}_{k_0
-1})\leq\rho^{1}_{n}(\mathbf{t}_{k_{0}-1})$, we obtain
\[
\rho_{n}^{2}(\mathbf{t}_{k_0})-\|c\|_{\infty}\leq\rho
_{n}^{1}(\mathbf{t}_{k_0}).
\]
Similarly, for $k<k_n$, $\rho_{n}^{m_{k}}(\mathbf{t}_{k})+\|c\|
_{\infty}
\geq\rho_{n}^{m_{k}}(\mathbf{t}_{k+1})>\rho_{n}^{m_{k}}(\mathbf
{t}_{k})$. Thus,
from
$\rho^{3}_{n}(\mathbf{t}_{k_{0}-1})\leq\rho^{2}_{n}(\mathbf
{t}_{k_0})=\rho^{2}
_{n}(\mathbf{t}_{k_{0}-1})$, we have
\[
\rho_{n}^{3}(\mathbf{t}_{k_0})\leq\|c\|_{\infty}+\rho
_{n}^{2}(\mathbf{t}_{k_0}).
\]
We have proved so far that the inequalities in (\ref{eq:recineq1})
hold for all $k\leq k_0$. Since
$|\rho^l_{n}(\mathbf{t}_{k+1})-\rho^l_{n}(\mathbf{t}_{k})|\leq\|c\|
_{\infty}$
and
$\rho^1_{n}(\mathbf{t}_{k_{0}-1})-\rho^3_{n}(\mathbf
{t}_{k_{0}-1})>2\|c\|_\infty$
we get
\[
\rho^1_{n}(\mathbf{t}_{k_0})-\rho^3_{n}(\mathbf{t}_{k_0})>0.
\]
Thus $m_{k_0}=3$ and $\rho_{n}^{3}(\mathbf{t}_{k_0})\leq\rho_{n}^{1}(
\mathbf{t}_{k_0})\leq\rho_{n}^{2}(\mathbf{t}_{k_0})$. Define
\[
k_{1}=\min\bigl\{k_n,\min\{k>k_{0}\dvtx m_{k}\neq3\mbox{ or } M_{k}\neq
2\}\bigr\}.
\]
For $k=k_{0},\ldots,k_{1}-1$,
$\rho^{2}_{n}(\mathbf{t}_{k+1})<\rho^2_n(\mathbf{t}_k)$ and $\rho
^1_n(\mathbf{t}_{k+1})
= \rho^1_n(\mathbf{t}_k)$ is constant, so the left-hand side inequality
of (\ref{eq:recineq1}) holds. Also, since $k_{1}\leq k_{n}$, for
$k\in\{k_{0}+1,\ldots,k_{1}-1\}$,
$\rho^{3}_{n}(\mathbf{t}_k)<\rho^{2}_{n}(\mathbf{t}_k)+4\|c\|
_{\infty}$. So
finally, (\ref{eq:recineq1}) holds for $k=0,\ldots,k_1$. Moreover,
if $k_{1}<k_{n}$, then $M_{k_1}=1$ and $m_{k_1}=3$. Indeed, as
above, $\rho_{n}^{2}(\mathbf{t}_{k_{1}-1})-\rho_{n}^{3}(\mathbf
{t}_{k_{1}-1})>2
\|c\|_{\infty}$ implies
\[
\rho_{n}^{2}(\mathbf{t}_{k_1})>\rho_{n}^{3}(\mathbf{t}_{k_1}).
\]
So $M_{k_1}\neq3$ and $m_{k_1}\neq2$. If $m_{k_1}=1$ and
$M_{k_1}=2$, then we write, by (\ref{eq:recineq1}),
\[
\rho_{n}^{1}(\mathbf{t}_{k_1})+\|c\|_{\infty}\geq\rho
_{n}^{2}(\mathbf{t}_{k_1})>
\rho_{n}^{3}(\mathbf{t}_{k_1})\geq\rho_{n}^{1}(\mathbf{t}_{k_1}).
\]
So $k_{1}=k_n$, a contradiction. Therefore, we necessarily have
$M_{k_1}=1$ and $m_{k_1}=3$. By recursion, it shows that for all
$k<k_n$, $m_{k}=3$. Hence, at each step one point is added to the
bin at $B_3$. No point is added to the bins at $B_1$ and $B_2$,
points may only be removed from the bins at $B_1$ and $B_2$. Since
there are at most $n$ points, we deduce $k_{n}\leq n$, as
claimed. Also, since $\Gamma^{l}(\mathbf{t}_0)=\mathbb{T}_l$, we obtain,
for all
$k=1, \ldots,k_n$, $\mathbb{T}_{3}\subset\Gamma^{3}(\mathbf{t}_k)$,
$\mathbb{T}_{2}\supseteq\Gamma^{2}(\mathbf{t}_k)$ and $\mathbb
{T}_{1}\supset
\Gamma^{1}(\mathbf{t}_k)$. The other case, where $m_{k_0}=2$ could be treated
similarly. So more generally, if, at some step, $l=m_k$ then $l
\neq M_j$ for all $k<j<k_n$, and conversely, if $l=M_k$ then $l
\neq m_j$ for all $k<j<k_n$. It implies that $\Gamma^{l}(\mathbf
{t}_k)$ is a
monotone sequence in $k$. Since $\Gamma^{l}(\mathbf{t}_0)=\mathbb{T}_l$,
for all $l
\in\{1,2,3\}$,
%
%
\begin{equation}\label{eq:monoD}
\Gamma^{l}_{n}\subseteq\mathbb{T}_{l} \quad\mbox{or}\quad \mathbb{T}
_{l}\subseteq
\Gamma^{l}_n.
\end{equation}
Assume now, that $t^1_{n}>zn^{-\alpha}$ with $z>0$ then, from
(\ref{eq:monoD}), $\mathbb{T}_1 \subseteq\Gamma^1_n$ and $\Gamma^2_n
\subseteq
\mathbb{T}_2$. For $t\in\mathbb{R}$, define the set
$V^{1}(t)=\Gamma^{1}(t,0,0)\setminus\mathbb{T}_1$. On the event
$\{t_{n}^{1}>zn^{-\alpha}\}$ we have
\[
\rho_{n}^{1}(\mathbf{t}_{n}) \geq n\int_{\mathbb{T}_1}c(x)\mu
_{n}(dx)+n
\int_{V^{1}(zn^{-\alpha})}c(x)\mu_{n}(dx)
\]
and
\[
\rho_{n}^{2}(\mathbf{t}_{n}) \leq
n\int_{\mathbb{T}_2}c_{2}(x)\mu_{n}(dx).
\]
So, by inequality (\ref{eq:hatAload}), we deduce that
on $\{t_{n}^{1}>zn^{-\alpha}\}$
\[
\int_{\mathbb{T}_1}c(x)\mu_{n}(dx)+\int_{V^{1}(zn^{-\alpha})}c
(x)\mu_{n}(dx)\leq\int_{\mathbb{T}_2}c_{2}(x)\mu_{n}(dx)
+\frac{2\|c\|_{\infty}}{n}.
\]
Or, equivalently,
%
%
\begin{eqnarray}\label{eq:Dx}\quad
\{t_{n}^{1}>zn^{-\alpha}\}& \subseteq& \biggl\{\sqrt{n}\int
_{V^{1}(zn^{-\alpha})}c(x)\mu_{n}(dx)\nonumber\\[-8pt]\\[-8pt]
&&\hspace*{5.7pt}\leq
\int_{\mathbb{T}_2}c_{2}(x)W_{n}(dx) -\int
_{\mathbb{T}
_1}c(x)W_{n}(dx)
+\frac{2\|c\|_{\infty}}{\sqrt{n}} \biggr\}.\nonumber
\end{eqnarray}
Let $A$ be a Borel set in $\mathbb{T}$. By Hoeffding's concentration
inequality (see, e.g., Corollary 2.4.14 in \cite{dembo}) we have,
for all $s\geq0$ and $l\in\{1,2,3\}$,
%
%
\begin{eqnarray}
\label{eq:mcdia}
P \biggl(\int_{A}c_{l}(x)\mu_{n}(dx)-\int_{A}c_{l}
(x) \,dx\geq s \biggr)&\leq&\exp(-K_{0}s^{2}n),\\
\label{eq:mcdiabis}
P \biggl(\int_{A}c_{l}(x)\mu_{n}(dx)-\int_{A}c_{l}
(x) \,dx\leq-s \biggr)&\leq&\exp(-K_{0}s^{2}n),
\end{eqnarray}
where $K_{0}=2\|c\|_{\infty}^{-2}$. Taking
$s=yn^{-\alpha}$, where $y>0$, we have
%
%
\begin{eqnarray}\label{eq:mcT1}
P \biggl(\int_{\mathbb{T}_l}c_{l}(x)W_{n}(dx)\geq y
n^{{1}/{2}-\alpha} \biggr)&\leq&\exp(-K_{0}y^{2}n^{1-2\alpha
}),\nonumber\\[-8pt]\\[-8pt]
P \biggl(\int_{\mathbb{T}_l}c_{l}(x)W_{n}(dx)\leq-yn^{{1}/
{2}-\alpha} \biggr)&\leq&\exp(-K_{0}y^{2}n^{1-2\alpha}).\nonumber
\end{eqnarray}
Similarly, by (\ref{eq:mcdiabis}) we deduce, for $s\geq0$,
\[
P\biggl(\int_{V^{1}(zn^{-\alpha})}c(x)\mu_{n}(dx)\leq
\int_{V^{1}(zn^{-\alpha})}c(x) \,dx-s \biggr)\leq\exp(
-K_{0}s^{2}n).
\]
By assumption (\ref{ass3}), there exists $c_{0}>0$ such
that $c(x)>c_{0}$, for all $x\in V^{1}(zn^{-\alpha})$. If $0\leq s
\leq\lambda\sqrt{3}/2$, the area of $V^{1}(s)$ is equal to
$\lambda s/4$. Therefore, for all $0\leq
z\leq\lambda\sqrt{3}n^{\alpha}/2$,
\[
K_{1}zn^{-\alpha}\leq\int_{V^{1}(zn^{-\alpha})}c(x) \,dx
\leq K_{2}zn^{-\alpha}
\]
with $K_{1}=c_{0}\lambda/4$ and $K_{2}=\|c\|_{\infty}
\lambda/4$. So, taking $s=K_{1}zn^{-\alpha}/2$, we get, for all $0
\leq z\leq\lambda\sqrt{3}n^{\alpha}$,
%
%
\begin{equation} \label{eq:mcV1}\qquad
P \biggl(\sqrt{n}\int_{V^{1}(zn^{-\alpha})}c(x)\mu_{n}(dx)
\leq\frac{K_{1}}{2}zn^{{1}/{2}-\alpha} \biggr)\leq
\exp(-K_{3}z^{2}n^{1-2\alpha}),
\end{equation}
where $K_{3}=K_{0}K_{1}^{2}/4$. Similarly, for $t\geq
0$, if $e_1 = ( 1, 0,0)$, $e_2 = (0,1,0)$, $e_3 = (0,0,1)$, we define
\[
U^{l}(t)= \bigl(\Gamma^{l}(t e_l)\setminus\mathbb{T}_{l} \bigr)\cup
\bigl(
\Gamma^{l}( - t e_l)\setminus\mathbb{T}_{\sigma(l)} \bigr),
\]
where $\sigma=(1\ 2\ 3)$ is the cyclic permutation. By
(\ref{eq:mcdia}) we have, for all $s\geq0$,
\[
P \biggl(\int_{U^{1}(zn^{-\alpha})}c(x)\mu_{n}(dx)\geq\int
_{U^{1}(zn^{-\alpha})}c
(x) \,dx+s \biggr)\leq\exp(-K_{0}s^{2}n).
\]
Thus, setting $s=zn^{-\alpha}$, we get
%
%
\begin{equation}\label{eq:mcV11}
P \bigl(\mu_{n}(U^{1}(zn^{-\alpha}))\geq K_4 zn^{-\alpha} \bigr)
\leq\exp(-K_{0}z^{2}n^{1-2\alpha})
\end{equation}
with $K_4 = 1 + 2 K_{2}$. Now, note that by
(\ref{eq:Dx}), from the union bound, for $y>0$,
\begin{eqnarray*}
\{t_{n}^{1}>zn^{-\alpha}\}
&\subseteq& \biggl\{\sqrt{n}\int_{V^{1}(zn^{-\alpha})}c(x)\mu
_{n}(dx)
\leq yn^{{1}/{2}-\alpha} \biggr\}\\
&&{}\cup
\biggl\{-\int_{\mathbb{T}_1}c_{1}(x)W_{n}(dx)+\frac{\|c\|_\infty
}{\sqrt
{n}}>\frac{1}{2}yn^{{1}/{2}-\alpha} \biggr\}
\\
&&{} \cup\biggl\{
\int_{\mathbb{T}_2}c_{2}(x)W_{n}(dx)+
\frac{\|c\|_\infty}{\sqrt{n}}>\frac{1}{2}yn^{{1}/{2}-\alpha
} \biggr\}.
\end{eqnarray*}
Now take $y=K_{1}z/2$. By (\ref{eq:mcT1}) and
(\ref{eq:mcV1}), if $4\|c\|_{\infty}n^{\alpha-1}{K_1}^{-1}\leq
z\leq\lambda\sqrt{3} n^{\alpha}$ we deduce
\begin{eqnarray*}
P(t^{1}_{n}>zn^{-\alpha})&\leq&\exp(-K_3
z^{2}n^{1-2\alpha} )+
2\exp\biggl(-\frac{K_{0}}{16}n^{1-2\alpha} (K_{1}z-4\|c\|
_{\infty}n
^{\alpha-1} )^{2} \biggr)\\
&\leq&3\exp\bigl(-K_{5}n^{1-2\alpha} (K_{1}z-4\|c\|_{\infty}
n ^{\alpha-1} )^{2} \bigr)
\end{eqnarray*}
with $K_{5}=\min\{K_{3}K_1^{-2},K_0/16\}$. Therefore, by
symmetry, for all $n$ and $z>0$ such that $4\|c\|_{\infty}n^{
\alpha-1}K_{1}^{-1}\leq z\leq\lambda\sqrt{3}n^{\alpha}/2$
%
%
\begin{equation}\label{eq:tailAhat}
P \Bigl(\max_{1\leq l\leq3}|t^l_{n}|>zn^{-\alpha} \Bigr)\leq18
e^{-K_{5}n^{1-2\alpha} (K_{1}z-4\|c\|_{\infty}n
^{\alpha-1} )^{2}}.
\end{equation}
Note that $\hat{\rho}_n=\rho_n(\hat{A})=\max_{1\leq
l\leq3}\rho_n^{l}(t_{n}^{l})$, so by (\ref{eq:hatAload}) we
have
\[
3\hat{\rho}_{n}-4\|c\|_{\infty}\leq\rho^1_n(\mathbf{t}_{n})+\rho^2_n
(\mathbf{t}_{n})+\rho^3_n(\mathbf{t}_{n})\leq3\hat{\rho}_n.
\]
Subtracting $3\sqrt{n}\gamma$, it follows
\[
3\frac{\hat{\rho}_{n}-n\gamma}{\sqrt
n}-\frac{4\|c\|_{\infty}}{\sqrt n}\leq\sqrt{n}\sum_{l=1}^{3}
\biggl(\int_{\Gamma^l_n}c_{l}(x)\mu_n(dx)-\gamma\biggr)
\leq
3\frac{\hat{\rho}_{n}-n\gamma}{\sqrt n}.
\]
Then we subtract the quantity
\[
\sum_{l=1}^{3}\int_{\mathbb{T}_l}c_{l}(x)W_{n}(dx)=\sqrt{n}
\sum_{l=1}^{3} \biggl(\int_{\mathbb{T}_l}c_{l}(x)\mu_{n}(dx)-\gamma
\biggr)
\]
and we get
%
%
\begin{eqnarray}\label{eq:hatrhon}\quad
&& \Biggl|3\frac{\hat{\rho}_n- n \gamma}{\sqrt{n}}-\sum_{l=1}^{3}
\int_{\mathbb{T}_l}c_{l}(x)W_{n}(dx) \Biggr| \nonumber\\[-8pt]\\[-8pt]
&&\qquad\leq\sqrt{n} \Biggl|
\sum_{l=1}^{3}\int_{\Gamma^l_n}c_{l}(x)\mu_n(dx)-\sum_{l=1}^3
\int_{\mathbb{T}_l}c_{l}(x)\mu_{n}(dx) \Biggr|+\frac{4\|c\|
_{\infty
}}{\sqrt
n}.\nonumber
\end{eqnarray}
Set $c_{\min}(x)=\min(c_1(x),c_2(x),c_3(x))$, and note that if
$x\in\mathbb{T}_l$ then $c_{\min}(x)=c_l(x)$. If $t^l_{n}\geq0$, we set
$V_{n}^{l}=V^{l}(t_{n}^{l})=\Gamma_{n}^{l}\setminus\mathbb{T}_l$,
and, if $t^l
_{n}<0$, we set $V^{l}_{n}=\Gamma_{n}^{\sigma(l)}\setminus\mathbb{T}_{l}$,
where $\sigma=(1\ 2\ 3)$ is the cyclic permutation. So
%
%
\begin{eqnarray}\label{eq:cmin}
&&\sum_{l=1}^{3}\int_{\Gamma^l_n}c_{l}(x)\mu_{n}(dx)-
\sum_{l=1}^{3}\int_{\mathbb{T}_l}c_{l}(x)\mu_{n}(dx)\nonumber\\
&&\qquad= \sum_{l=1}^{3}\int_{\Gamma^l_n} \bigl(c_{l}(x)-c_{\min}(x)
\bigr)\mu
_{n}(dx)\nonumber\\[-8pt]\\[-8pt]
&&\qquad=\sum_{l=1}^{3}\mathbh{1}\{t_{n}^{l}\geq0\}
\int_{V^{l}_n} \bigl(c_{l}(x)-c_{\min}(x) \bigr)
\mu_{n}(dx) \nonumber\\
&&\qquad\quad{} +\sum_{l=1}^{3}\mathbh{1}\{t_{n}^{l}<0\}
\int_{V^l_n} \bigl(c_{\sigma(l)}(x)-c_{\min}(x) \bigr)
\mu_{n}(dx).\nonumber
\end{eqnarray}
Note that if $x\in\mathbb{T}_m$, with $m\neq l$, then $|c_{l}
(x)-c_{\min}(x)|=|c_{l}(x)-c_{m}(x)|$. For example, assume $l=1$,
$m=2$ and $x=t\rme^{i{\pi}/{6}+i\theta}\in\mathbb{T}_2$, with
$0\leq\theta\leq\pi/3$, we then have
\begin{eqnarray*}
|c_1(x)-c_{\min}(x)|&=&|c_{1}(x)-c_{2}(x)|=|c(t\rme^{i{\pi}/
{6}+i\theta})-c(t\rme^{i{\pi}/{6}+i\theta}\rme^{-i
{2\pi}/{3}})|\\
&=&|c(t\rme^{i{\pi}/{6}+i\theta})-c(t\rme^{-i
{\pi}/{2}+i\theta})|.
\end{eqnarray*}
By the symmetry assumption (\ref{ass4}), we deduce
\[
|c_{1}(x)-c_{\min}(x)|=|c(t\rme^{i{\pi}/{6}
+i\theta})-c(t\rme^{i{\pi}/{6}-i\theta})|.
\]
Again by assumption (\ref{ass4}), $c$ is Lipschitz in a
neighborhood of $D_{1}\cup D_{3}$. Letting $L>0$ denote the
Lipschitz constant, if $x$ is close enough to $D_1$, say the
distance $d(x,D_1)$ from $x$ to $D_1$ is less than or
equal to $\varepsilon$ with $0<\varepsilon<\lambda\sqrt{3}/2$, we
have
\begin{eqnarray*}
|c_1(x)-c_{\min}(x)|&\leq& Lt|\rme^{i{\pi}/{6}
+i\theta}-\rme^{i{\pi}/{6}-i\theta}|
=Lt|\rme^{i\theta}-\rme^{-i\theta}|\\
&=&2Lt\sin\theta=2Ld(x,D_1).
\end{eqnarray*}
By symmetry, for all $l\in\{1,2,3\}$, if
$d(x,D_l)\leq\varepsilon$, then
\[
|c_{l}(x)-c_{\min}(x)|\leq2Ld(x,D_l) \quad\mbox{and}\quad
\bigl|c_{\sigma(l)}(x)-c_{\min}(x)\bigr|\leq2Ld(x,D_l).
\]
Fix $\alpha\in(1/4,1/2)$, $z>0$ and choose $n$ large
enough so that $4\|c\|_{\infty}n^{\alpha-1}K_{1}^{-1}\leq z\leq
\varepsilon n^{\alpha}$. Then, by (\ref{eq:tailAhat}) with
probability at least $1-18e^{-K_{5}n^{1-2\alpha} (K_{1}
z-4\|c\|_{\infty}n^{\alpha-1} )^2 }$, we have
$\max_{1\leq l\leq3}|t^l_{n}|\leq zn^{-\alpha}$. On this event,
if $x\in V^l(t_{n}^{l})$ then $d(x,D_l)\leq zn^{-\alpha}
\leq\varepsilon$. It follows by (\ref{eq:cmin}) that, with
probability at least $1-18\times\break e^{-K_{5}n^{1-2\alpha} (K_{1}
z-4\|c\|_{\infty}n^{\alpha-1} )^{2}}$,
\begin{eqnarray*}
&&
\sqrt{n} \Biggl|\sum_{l=1}^{3}\int_{\Gamma^l_n}c_{l}(x)\mu
_{n}(dx)-
\sum_{l=1}^{3}\int_{\mathbb{T}_l}c_{l}(x)\mu_{n}(dx)
\Biggr|\\
&&\qquad\leq
\sqrt{n}\sum_{l=1}^{3}2Lzn^{-\alpha}\mu_{n}(V^{l}_n)\\
&&\qquad\leq2Lzn^{{1}/{2}-\alpha}\sum_{l=1}^{3}\mu_{n}(U^{l}(zn
^{-\alpha})).
\end{eqnarray*}
By (\ref{eq:mcV11}), with probability at least $1-3
\exp(-K_{0}z^{2}n^{1-2\alpha})$, it holds\break $\sum_{l=1}^{3}
\mu_{n}(U^{l}(zn^{-\alpha}))\leq3K_{4}zn^{-\alpha}$. Using that
for all events $A,B$ it holds $P(A\cap B)\geq1-P(A^c)-P(B^c)$, we
obtain, for all $n$ large enough so that $4\|c\|_{\infty}n^{\alpha
-1}\times\break K_{1}^{-1}\leq z\leq\varepsilon n^{\alpha}$,
\[
\sqrt{n} \Biggl|\sum_{l=1}^{3}\int_{\Gamma^l_n}c_{l}(x)\mu
_{n}(dx)
-\sum_{l=1}^{3}\int_{\mathbb{T}_l}c_{l}(x)\mu_{n}(dx) \Biggr|
\leq12 LK_{4}z^{2}n^{{1}/{2}-2\alpha}
\]
with probability at least $1-21\exp(-K_{6}n^{1-2\alpha}
(K_{1}z-4\|c\|_{\infty}n^{\alpha-1} )^{2} )$, where
$K_{6}=\min\{K_{0}{K_1}^{-2},K_5\}$. By this latter inequality and
(\ref{eq:hatrhon}), with the same probability,
\[
\Biggl|3\frac{\hat{\rho}_{n}-\gamma}{\sqrt{n}}-\sum_{l=1}^{3}
\int_{\mathbb{T}_l}c_{l}(x)W_{n}(dx) \Biggr|\leq
12LK_{2}z^{2}n^{{1}/{2}-2\alpha}+4\|c\|_{\infty}n^{-{1}/{2}}.
\]
Fix $z=(24LK_2)^{-{1}/{2}}$ so that $12LK_{2}z^{2}
=1/2$. Then there exists $n_0$ such that, for all $n\geq n_0$, $4
\|c\|_{\infty}n^{\alpha-1}K_{1}^{-1}\leq z\leq\varepsilon n
^{\alpha}$ and $8\|c\|_{\infty}n^{-{1}/{2}}\leq
n^{{1}/{2}-2\alpha}$. Then, for all $n\geq n_0$,
%
%
\begin{equation}\label{eq:final31}
\Biggl|3\frac{\hat{\rho}_{n}-\gamma}{\sqrt{n}}-\sum_{l=1}^{3}
\int_{\mathbb{T}_l}c_{l}(x)W_{n}(dx) \Biggr|\leq n^{1/2-2\alpha}
\end{equation}
with probability at least
\begin{eqnarray*}
&&
1-21\exp\bigl(-K_{6}n^{1-2\alpha} \bigl(K_{1}(24LK_2)^{-{1}/{
2}}-4\|c\|_{\infty}n_{0}^{\alpha-1} \bigr)^{2} \bigr)\\
&&\qquad=1-K_{7}
\exp(-K_{8}n^{1-2\alpha} ).
\end{eqnarray*}
Finally, we set $L_{0}=K_8$ and $L_{1}=\max\{K_7,K_9\}$, where
$K_9=\exp(K_{8}n_0^{1-2\alpha} )$. With this choice of
$L_0$ and $L_1$, (\ref{eq:final31}) holds for all $n\geq1$ with
probability at least $1-L_{1}\exp(-L_{0}n^{1-2\alpha})$.
\end{pf*}

\section{Large deviation principles}\label{sec:LDP}

In this section we provide LDPs for the optimal and suboptimal
load. Letting $\ll$ denote absolute continuity between measures,
we define by
\[
H(\nu| \ell)=\cases{
\displaystyle\int_{\mathbb{T}}\frac{d\nu}{d\ell}(x)\log\frac
{d\nu
}{d\ell}(x)\,
d\ell, &\quad if $\nu\ll\ell$,\vspace*{2pt}\cr
+\infty, &\quad otherwise,}
\]
the relative entropy of $\nu\in\mathcal{M}_{1}(\mathbb{T})$
with respect to the Lebesgue measure $\ell$. Moreover, if $f$ is a
nonnegative measurable function on $\mathbb{T}$, we denote by $\ell_f$
the measure on $\mathbb{T}$ with density $f$. In particular, if
$\int_{\mathbb{T}}f(x) \,dx=1$, we set
\[
H(f)=H(\ell_f | \ell)=\int_\mathbb{T}f(x)\log f(x) \,dx.
\]

\subsection{Combining Sanov's theorem and the contraction principle}

Next Theorem \ref{th:LDPc} follows combining Sanov's theorem and
the contraction principle.
\begin{Theorem}\label{th:LDPc}
Assume (\ref{ass20}) and (\ref{ass2}). Then:

\begin{longlist}
\item $\{\rho_n/n\}_{n\geq1}$ satisfies an LDP on
$\mathbb{R}$ with good rate function
%
%
\begin{equation}\label{eq:rf}
J(y)=\inf_{\alpha\in\mathcal{M}_1(\mathbb{T})\dvtx\Phi(\alpha
)=y}H(\alpha|
\ell).
\end{equation}

\item $\{\overline{\rho}_n/n\}_{n\geq1}$ satisfies an
LDP on $\mathbb{R}$ with good rate function
%
%
\begin{equation}\label{eq:ratef}
\overline{J}(y)=\inf_{\alpha\in\mathcal{M}_1(\mathbb{T})\dvtx
\Psi(\alpha
)=y}H(\alpha|
\ell).
\end{equation}
\end{longlist}
\end{Theorem}
\begin{pf}
By Sanov's theorem (see, e.g., Theorem 6.2.10 in \cite{dembo})
the sequence $\{\frac{1}{n}\sum_{i=1}^{n}\delta_{X_i}\}_{n\geq1}$
satisfies an LDP on $\mathcal{M}_1(\mathbb{T})$, with good rate function
$H(\cdot| \ell)$. Recall that the space $\mathcal{M}_{1}(\mathbb{T})$,
equipped with the topology of weak convergence, is a Hausdorff
topological space (refer to \cite{billingsley}). By Lemma
\ref{le:tvd} the function $\Phi$ is continuous on
$\mathcal{M}_{1}(\mathbb{T})$. Therefore, using (\ref{eq:rhophi})
and the
contraction principle (see, e.g., Theorem 4.2.1 in \cite{dembo}) we
deduce that the sequence $\{\widetilde{\rho}_n/n\}_{n\geq1}$
satisfies an LDP on $\mathbb{R}$ with good rate function given by
(\ref{eq:rf}). Consequently, by Lemma \ref{le:tilderho}(iii)
and Theorem 4.2.13 in \cite{dembo}, $\{\rho_n/n\}_{n\geq1}$ obeys
the same LDP. The proof of (ii) is identical and follows from
(\ref{eq:phirho}).
\end{pf}
\begin{Remark}
It is worthwhile noticing that one can prove Theorem
\ref{th:LDPc} also by applying Lemmas \ref{le:tvd},
\ref{le:tilderho}\textup{(iii)} and the results in O'Connell
\cite{oconnell}.
\end{Remark}

\subsection{Computing $\Lambda^*$ and $\overline{\Lambda}{}^*$}
\label{subsec:laplace}

In this subsection we compute the Fenchel--Legendre transforms
$\Lambda^*$ and $\overline{\Lambda}{}^*$.

\subsubsection{\texorpdfstring{Proof of Proposition \protect\ref{prop:Laplace}}{Proof of Proposition 1.4}}

We only compute $\Lambda^{*}$ in (i). The expression of
$\overline{\Lambda}{}^{*}$ in (ii) can be computed similarly.
Clearly, for $\theta\in\mathbb{R}$,
\[
\Lambda'(\theta)=\frac{\int_{\mathbb{T}_1}c(x)\rme^{\theta
c(x)} \,dx}{\int_{\mathbb{T}_1} \rme^{\theta
c(x)} \,dx}
\]
and
\[
\Lambda''(\theta)=\int_{\mathbb{T}_1} c^2(x)\frac{\rme^{\theta
c(x)}}{\int_{\mathbb{T}_1}\rme^{\theta
c(x)} \,dx} \,dx- \biggl(\int_{\mathbb{T}_1}c(x)\frac
{\rme^{\theta
c(x)}}{\int_{\mathbb{T}_1}\rme^{\theta
c(x)} \,dx} \,dx \biggr)^2>0
\]
[the strict inequality comes from the assumption that $c(\cdot)$ is not
constant on $\mathbb{T}_1$]. Therefore, the function $\Lambda'$ is
strictly increasing. Consider the probability measure on $\mathbb{T}_1$:
\[
\mathrm{P}_\theta(dx)=\frac{\rme^{\theta
c(x)} \,dx}{\int_{\mathbb{T}_1}\rme^{\theta
c(x)} \,dx}.
\]
Next Lemma \ref{le:Laplace} is classical; we give a
proof for completeness.
\begin{Lemma}\label{le:Laplace}
Under the assumptions of Proposition \ref{prop:Laplace}, the
following weak convergence holds:
\[
\mathrm{P}_{\theta}\Rightarrow\delta_{0} \qquad\mbox{as }
\theta\to+\infty\quad\mbox{and}\quad \mathrm{P}_\theta\Rightarrow
\delta
_{B_1} \qquad\mbox{as }
\theta\to-\infty.
\]
\end{Lemma}
\begin{pf}
We only prove the first limit. Indeed, the
second limit can be showed similarly. We need to show
\[
\mathrm{P}_\theta(A)\to\delta_0(A) \qquad\mbox{as }\theta\to+\infty
\mbox{ for
any Borel set }A\subseteq\mathbb{T}_1\mbox{ such that }0\notin
\partial A.
\]
If $0\notin A\subseteq\mathbb{T}_1$ then, by assumption
(\ref{eq:concave}), $c(x)<c(0)$ for any $x\in A$. So $A\subseteq
I_t$, for some $t>0$, where $I_{t}=\{x\in\mathbb{T}_1\dvtx c(x)\leq
c(0)-t\}$. By assumption $c$ is continuous at $0$, so there exists
an open neighborhood of $0$, say $V_t$, such that, for all $x\in
V_t$, $c(x)\geq c(0)-t/2$. Note that, for any $\theta>0$,
\begin{eqnarray*}
\mathrm{P}_\theta(I_t)&=&\int_{I_t}\frac{\rme^{\theta
c(x)}}{\int_{\mathbb{T}_1}\rme^{\theta
c(x)} \,dx} \,dx\\
&\leq&\int_{\mathbb{T}_1}\frac{\rme^{\theta c(0)-\theta t
}}{\int_{V_t\cap\mathbb{T}_1}\rme^{\theta c(0)-\theta t
/2} \,dx} \,dx\\
&\leq&\ell(V_t\cap\mathbb{T}_1)^{-1}\rme^{-\theta t /2}.
\end{eqnarray*}
Thus, for all $t>0$, $\lim_{\theta\to+\infty}\mathrm{P}_\theta
(I_t)=0$. This guarantees the claim in the case when the Borel set
$A\subseteq\mathbb{T}_1$ does not contain $0$. Suppose now $0\in A$, then
$0 \notin\mathbb{T}_1 \setminus A$, and we get $\mathrm{P}_\theta
(A) = 1 -
\mathrm{P}_\theta( \mathbb{T}_1 \setminus A ) \to1$ as $\theta$
goes to
infinity.
\end{pf}

We can now continue the proof of the proposition. Let $c(B_1)<y<
c(0)$. By Lemma 2.3.9(b) in \cite{dembo}, we need to show that
there exists a unique solution $\theta_y$ of $\Lambda'(\theta)=
y$. To this end, note that
$\Lambda'(\theta)=\int_{\mathbb{T}_1}c(x) \mathrm{P}_{\theta
}(dx)$. By
assumption, $c$ is continuous at $0$ and $B_1$, so by Lemma
\ref{le:Laplace} and Theorem 5.2 in \cite{billingsley} it follows
\[
\lim_{\theta\to-\infty}\Lambda'(\theta)=
c(B_1)<y<c(0)=\lim_{\theta\to+\infty}\Lambda'(\theta).
\]
Since $\Lambda'$ is continuous and strictly increasing, the mean
value theorem implies the existence and uniqueness of
$\theta_y$. Consider now $y>c(0)$. Note that, for $\theta\geq0$,
$\Lambda(\theta)\leq\theta c(0)$. Therefore
\[
\theta y-\Lambda(\theta) \geq\theta\bigl(y-c(0)\bigr).
\]
It follows that $\Lambda^*(y)=+\infty$. Similarly, for
$y<c(B_1)$, we use that, for $\theta\leq0$, $\Lambda(\theta)\leq
\theta c(B_1)$ and deduce $\Lambda^*(y)=+\infty$. Finally we prove
(iii). We first show that
%
%
\begin{equation}\label{eq:inloglap}
\Lambda(\theta/3)<\overline{\Lambda}(\theta)\qquad \mbox{for all }
\theta>0.
\end{equation}
Showing (\ref{eq:inloglap}) amounts to show that, for all
$\theta>0$,
%
%
\begin{equation}\label{ineq:equiv}
\int_{\mathbb{T}_1}\rme^{\theta c(x)} \,dx + 2/3 - 3
\int_{\mathbb{T}_1} \rme^{\theta c(x)/3} \,dx > 0.
\end{equation}
By Jensen's inequality it follows that
\[
\biggl(\int_{\mathbb{T}_1}\rme^{\theta c(x)/3} \,dx \biggr)^3
<\frac{1}{9}\int_{\mathbb{T}_1}\rme^{\theta c(x)} \,dx
\]
(the strict inequality derives from the strict convexity
of the cubic power on $[0,\infty)$, and the fact that $c$ is not
constant on $\mathbb{T}_1$). Hence the left-hand side of
(\ref{ineq:equiv}) is larger than $9(\int_{\mathbb{T}_1}\rme
^{\theta c(x)/3} \,dx
)^{3}-3\int_{\mathbb{T}_1}\rme^{\theta c(x)/3} \,dx + 2/3$,
which is equal to
\[
9 \biggl(\int_{\mathbb{T}_1}\rme^{\theta
c(x)/3} \,dx-\frac{1}{3}\biggr)^{2}\biggl(
\int_{\mathbb{T}_1}\rme^{\theta c(x)/3} \,dx+\frac2 3
\biggr),
\]
and inequality (\ref{ineq:equiv}) follows. Now, let
$\gamma<y<c(0)/3$. By Theorem \ref{th:lln},\break $
\lim_{n\to\infty}\rho_n/n=\lim_{n\to\infty}\overline{\rho}_n/n
=\gamma<y$. Thus, by Lemma 2.2.5 in \cite{dembo} we have
\[
\Lambda^*(3y)=\sup_{\theta>0}\bigl(\theta y-\Lambda(\theta/3)\bigr)
\quad\mbox{and}\quad \overline{\Lambda}{}^*(y)=\sup_{\theta>0}\bigl(\theta
y-\overline{\Lambda}(\theta)\bigr)=\eta_y y-\overline{\Lambda}(\eta_y),
\]
where $\eta_y$ is the unique positive solution of
(\ref{eq:etay}). Finally, (\ref{eq:inloglap}) yields
\[
\overline{\Lambda}{}^*(y)=y\eta_y-\overline{\Lambda}(\eta_y)
<y\eta_y-\Lambda(\eta_y/3)\leq\sup_{\theta>0}\bigl(\theta
y-\Lambda(\theta/3)\bigr)=\Lambda^*(3y).
\]

\subsubsection{Value of the Fenchel--Legendre transforms at the
extrema}\label{subsub:laplaceextrema}

In this paragraph, for the sake of completeness, we deal with the
value of $\Lambda^*$ and $\overline{\Lambda}{}^*$ at $c(B_1)$ and
$c(0)$. If $c$ is differentiable as a function from $\mathbb{T}\subset
\mathbb{C}$
to $\mathbb{R}$, we denote by $\operatorname{grad}_x (c)$ its
gradient at $x$. The
following proposition holds:
\begin{Proposition}\label{prop:extrema}
Suppose that the assumptions of Proposition \ref{prop:Laplace}
hold and that $c$ is differentiable at $0$ and $B_1$. If, moreover,
for all $\omega\in[-\pi/2,\pi/6]$, $\operatorname{grad}_0
(c)\cdot\rme^{i\omega}<0$ and, for all
$\omega\in[2\pi/3,\pi]$, $\operatorname{grad}_{B_1}(c)\cdot\rme
^{i\omega}
>0$, then
\[
\Lambda^{*}(c(B_1))=\overline{\Lambda}{}^{*}(c(B_1))=\Lambda^*
(c(0))=\overline{\Lambda}{}^{*}(c(0))=+\infty.
\]
\end{Proposition}
\begin{pf}
We show the proposition only for $\Lambda^*(c(0))$. The other
three cases can be proved similarly. Using polar coordinates, we
have
\[
\int_{\mathbb{T}_1}\rme^{\theta c(x)} \,dx=\int_{-\pi/2}
^{\pi/6}\int_{I_\omega}\rme^{\theta
c(r\rme^{i\omega})}r \,dr\,d\omega
\]
for some segment $I_\omega=[0,a_\omega]$. Laplace's method
(see, e.g., Murray \cite{murray}) gives, for all
$\omega\in[-\pi/2,\pi/6]$,
\[
\int_{I_\omega}\rme^{\theta
c(r\rme^{i\omega})}r \,dr
\sim\frac{\rme^{\theta c(0)}}{\theta^2|{\operatorname{grad}_0}
(c)\cdot\rme^{i\omega}|} \qquad\mbox{as }
\theta\rightarrow+\infty,
\]
where we write $f\sim g$ if $f$ and $g$ are two
functions such that, as $x\to+\infty$, the ratio $f(x)/g(x)$
converges to $1$. We deduce that, as $\theta\rightarrow+\infty$,
\[
\int_{\mathbb{T}_1}\rme^{\theta c(x)} \,dx\sim
\rme^{\theta
c(0)}\theta^{-2}\int_{-\pi/2}^{\pi/6}\frac{1}{|{\operatorname{grad}_0}
(c)\cdot\rme^{i\omega}|} \,d\omega.
\]
Since the integral in the right-hand side is a finite
positive constant, we have $\Lambda(\theta)=\theta
c(0)-2\log\theta+o(\log\theta)$, and therefore
\[
\Lambda^*(c(0))=\sup_{\theta\in\mathbb{R}}\bigl(\theta
c(0)-\Lambda(\theta)\bigr)=\sup_{\theta\in\mathbb{R}}\bigl(2\log\theta
+o(\log
\theta)\bigr)=+\infty.
\]
\upqed\end{pf}

In the next two subsections, we solve some variational
problems. We refer the reader to the book by Buttazzo, Giaquinta
and Hildebrandt \cite{buttazzo} for a survey on calculus of
variations.

\subsection{\texorpdfstring{Proof of Theorem
\protect\ref{th:LDP}\textup{(i)}}{Proof of Theorem 1.3(i)}}

We divide the proof of Theorem \ref{th:LDP}(i) in $5$ steps.

\subsubsection*{Step 1: Case $y\protect\notin(c(B_1)/3,c(0)/3)$}

We have to prove that $J(y)=\infty$. Denote by
$\mathcal{M}_1^{\mathrm{ac}}(\mathbb{T})\subseteq\mathcal{M}_1(\mathbb{T})$
the set of
probability measures on $\mathbb{T}$ which are absolutely continuous with
respect to $\ell$. For $\alpha\in\mathcal{M}_1^{\mathrm{ac}}(\mathbb{T})$, define
the measures in $\mathcal{M}_{b}(\mathbb{T})$
\[
\alpha_l(\,dx)=\frac{c_{\sigma^2(l)}(x)c_{\sigma
(l)}(x)}{c_1(x)c_2(x)+c_1(x)c_3(x)+c_2(x)c_3(x)}
\alpha(\,dx),\qquad l\in\{1,2,3\},
\]
where $\sigma=(1\ 2\ 3)$ is the cyclic
permutation. Clearly $\alpha_1+\alpha_2+\alpha_3=\alpha$ and
%
%
\begin{equation}\label{eq:abcin}
\Phi(\alpha)\leq\phi(\alpha_1,\alpha_2,\alpha_3)<c(0)/3,
\end{equation}
where the strict inequality follows by assumption
(\ref{eq:concave1}) and the fact that $\alpha$ is a probability
measure on $\mathbb{T}$ such that $\alpha\ll\ell$. The above argument
shows that
$\{\alpha\in\mathcal{M}_1^{\mathrm{ac}}(\mathbb{T})\dvtx\Phi(\alpha)=y\}
=\varnothing$,
for all $y\geq c(0)/3$. Therefore, by Theorem \ref{th:LDPc}(i),
we have $J(y)=+\infty$ if $y\geq c(0)/3$. Using assumptions
(\ref{ass3}) and (\ref{eq:concave}), one can easily realize
that, for any measure $\beta\in\mathcal{M}_b(\mathbb{T})$, $\beta
(c_l)\geq
c(B_1)\beta(\mathbb{T})$ and the equality holds only if
$\beta=\delta_{B_l}$. By Lemma \ref{le:tilderho}(i) we deduce
that, for all $\alpha\in\mathcal{M}_1(\mathbb{T})$, $3\Phi(\alpha
)>c(B_1)$. This
gives $J(y)=\infty$ for all $y\leq c(B_1)/3$, and concludes the
proof of this step.

\subsubsection*{Step 2: The set function $\nu$ and an alternative
expression for $\Lambda^*(3y)$}

For the remainder of the proof we fix $y\in(c(B_1)/3,
c(0)/3)$. For this we shall often omit the dependence on $y$ of
the quantities under consideration. In this step we give an
alternative expression for $\Lambda^*(3y)$ that will be used later
on. Let $B\subset\mathbb{T}$ be a Borel set with positive Lebesgue
measure. Define the function of $(\eta_0,\eta_1)\in\mathbb{R}^2$
\[
m(B,\eta_0,\eta_1)=\int_{B}\rme^{-1-\eta_0-\eta_1
c(x)} \,dx.
\]
It turns out that $m(B,\cdot)$ is strictly convex on
$\mathbb{R}^2$ (the second derivatives with respect to $\eta_0$ and
$\eta_1$ are strictly bigger than zero). Define the strictly
concave function
\[
F(B,\eta_0,\eta_1)=-\eta_0-3y\eta_1-3 m(B,\eta_0,\eta_1)
\]
and the set function
\[
\nu(B)=\sup_{(\eta_0,\eta_1)\in\mathbb{R}^2}F(B,\eta_0,\eta_1).
\]
Arguing as in the proof of Lemma 2.2.31(b) in
\cite{dembo}, we have
\begin{eqnarray*}
&&\operatorname{grad}_{(\gamma_0,\gamma_1)}(3m(B,\cdot))=(-1,-3y)\\
&&\quad
\Rightarrow\nu(B)=(\gamma_0,\gamma_1)\cdot(-1,-3y)-3
m(B,\gamma_0,\gamma_1),
\end{eqnarray*}
where $\cdot$ denotes the scalar product on
$\mathbb{R}^2$. Therefore, if there exist $\gamma_0=\gamma_{0}(B)$ and
$\gamma_1=\gamma_{1}(B)$ such that
%
%
\begin{equation}\label{eq:gammaeq0}
\int_{B}\rme^{-\gamma_{1}
c(x)}\,dx=\rme^{1+\gamma_{0}}/3 \quad\mbox{and}\quad
\int_{B}c(x)\rme^{-\gamma_{1}
c(x)}\,dx=y\rme^{1+\gamma_{0}},
\end{equation}
then it is easily seen that
\[
\nu(B)=-\bigl(1+\gamma_0(B)\bigr)-3y\gamma_1(B).
\]
In particular, by Proposition \ref{prop:Laplace}(i),
setting $\gamma_1(\mathbb{T}_1)=-\theta_{3y}$ and
$\gamma_0(\mathbb{T}_1)=\Lambda(\theta_{3y})-1$, one has
%
%
\begin{equation}\label{eq:nuT1}
\Lambda^*(3y)=\nu(\mathbb{T}_1)=
-\bigl(1+\gamma_{0}(\mathbb{T}_1)\bigr)-3y\gamma_{1}(\mathbb{T}_1),
\end{equation}
and $\gamma_{0}(\mathbb{T}_1)$ and $\gamma_{1}(\mathbb{T}_1)$ are the
unique solutions of the equations in (\ref{eq:gammaeq0}) with
$B=\mathbb{T}_1$. Note also that, for Borel sets $A$ and $B$ such that
$A\subseteq B\subseteq\mathbb{T}$, we have for all $\eta_0,\eta
_1\in\mathbb{R}$,
\[
m(B,\eta_0,\eta_1)-m(A,\eta_0,\eta_1)=\int_{\mathbb{T}}\bigl(\mathbh{1}
_{B}(x)-\mathbh{1}_{A}(x)\bigr)
\rme^{-1-\eta_0-\eta_1 c(x)} \,dx\geq0.
\]
In particular, for all $\eta_0,\eta_1\in\mathbb{R}$,
$F(A,\eta_0,\eta_1)\geq F(B,\eta_0,\eta_1)$. This proves that the
set function $\nu$ is nonincreasing (for the set inclusion). An
easy consequence is the following lemma. For $B\subset\mathbb{T}$ and
$z\in\mathbb{C}$, define $zB=\{zx\dvtx x\in B \}$ and
\begin{eqnarray*}
\mathcal{T}&=&\bigl\{\mbox{Borel sets }
B\subset\mathbb{T}\dvtx\ell(B)>0 \mbox{ and}\\
&&\hspace*{5.1pt} \ell\bigl(B\cap(jB)\bigr)=\ell
\bigl(B\cap
(j^{2}B)\bigr)=\ell\bigl((jB)\cap(j^{2}B)\bigr)=0\bigr\}.
\end{eqnarray*}

\begin{Lemma}\label{le:minmax}
Under the foregoing assumptions and notation, it holds
\[
\inf\{\nu(B)\dvtx B\in\mathcal{T}\}=\inf\{\nu(B)\dvtx B\in
\mathcal{T}\mbox{ and }
\ell(B)=1/3\}< +\infty.
\]
\end{Lemma}
\begin{pf}
The monotonicity of $\nu$ implies
$\nu(\mathbb{T})\leq\nu(\mathbb{T}_1)$. So the finiteness of the
infimum follows
by $\nu(\mathbb{T}_1)<+\infty$ that we proved above. Note that if
$B\in
\mathcal{T}$, then $B\cup(jB)\cup(j^{2}B)\subset\mathbb{T}$ and
$1\geq
\ell(B\cup(jB)\cup(j^{2}B))=\ell(B)+\ell(jB)+\ell(j^2 B)=3
\ell(B)$. So
\[
\inf\{\nu(B)\dvtx B\in\mathcal{T}\}=\inf\{\nu(B)\dvtx B\in
\mathcal{T}\mbox{ and }
\ell(B)\leq1/3\}.
\]
Now, if $B\in\mathcal{T}$ is such that $\ell(B)<1/3$, define the
set $C=\mathbb{T}\setminus(B\cup(jB)\cup(j^{2}B))$; note that
$\ell(C)=1-3\ell(B)>0$ and $C=jC=j^2 C$. Set $C_1=C\cap\mathbb
{T}_1$ and
define $D=B\cup C_1$. Clearly, $B\subset D$ and therefore
$\nu(B)\geq\nu(D)$. Moreover, it is easily checked that
$D\in\mathcal{T}$. Indeed, $\ell(D)\geq\ell(B)>0$ and, for instance,
\begin{eqnarray*}
\ell\bigl(D\cap(jD)\bigr)&=&\ell\bigl((B\cup C_1)\cap\bigl((jB)\cup(jC_1)\bigr)\bigr)\\
&\leq&\ell\bigl(B\cap(jB)\bigr)+\ell\bigl(B\cap
(jC_1)\bigr)+\ell\bigl(C_1\cap(jB)\bigr)+\ell\bigl(C_1\cap(jC_1)\bigr)\\
&=&0.
\end{eqnarray*}
The claim follows since
\[
\ell(D)=\ell(B)+\ell(C_1)=\ell(B)+\ell(C)/3=1/3.
\]
\upqed\end{pf}

\subsubsection*{Step 3: The related variational problem}

As above, we fix $y\in(c(B_1)/3$, $c(0)/3)$. Recall that
$H(\alpha| \ell)=+\infty$ if $\alpha$ is not absolutely
continuous with respect to $\ell$. So, by Theorem
\ref{th:LDPc}(i),
\[
J(y)=\inf_{\alpha\in\mathcal{M}_1^{\mathrm{ac}}(\mathbb{T})\dvtx\Phi
(\alpha)=y}
H(\alpha| \ell).
\]
Define the following functional spaces:
\[
\mathcal{B}=\{\mbox{measurable functions defined on }\mathbb{T}\mbox
{ with values in }
[0,\infty)\}
\]
and
\begin{eqnarray*}
\mathcal{B}^3_{\Phi}&=& \Biggl\{(f_1,f_2,f_3)\in\mathcal{B}^{3}\dvtx\ell
\Biggl(\sum
_{l=1}^{3}f_{l} \Biggr)=1\mbox{
and}\\
&&\hspace*{7.9pt}\phi(\ell_{f_1},\ell_{f_2},\ell_{f_3})=\Phi(\ell_{f_1}
+\ell_{f_2}+\ell_{f_3}) \Biggr\}
\end{eqnarray*}
(recall that $\ell_f$ is the measure with density $f$). By Lemma
\ref{le:tilderho}(i) it follows
%
%
\begin{equation}
\label{eq:varpb} J(y)=\inf_{(f_1,f_2,f_3) \in\mathcal{R}_\Phi^3}
H \Biggl(\sum_{l=1}^3 f_l (x) \Biggr),
\end{equation}
where
\[
\mathcal{R}^3_{\Phi}= \{(f_1,f_2,f_3)\in\mathcal{B}^3_{\Phi}\dvtx
\phi(
\ell_{f_1},\ell_{f_2},\ell_{f_3})=y \}
\]
(note that the superscript ``3'' in $\mathcal{B}^{3}_{\Phi}$ and $
\mathcal{R}^{3}_{\Phi}$ is a reminder that these spaces are defined on
triplets of functions in $\mathcal{B}$; it is not related to the Cartesian
product of three spaces). Computing the value of $J(y)$ from
(\ref{eq:varpb}) is far from obvious; indeed $\mathcal{R}^3_{\Phi}$
is not
a convex set, and the standard machinery of calculus of variations
cannot be applied directly. The key idea is the following:
consider the same minimization problem on a larger convex space,
defined by linear constraints; compute the solution of this
simplified variational problem; show that this solution is in
$\mathcal{R}^3_{\Phi}$. To this end, note that, again by Lemma
\ref{le:tilderho}(i), if $(f_1,f_2,f_3)\in\mathcal{B}^3_{\Phi}$, then
$\ell_{f_1}(c_1)=\ell_{f_2}(c_2)=\ell_{f_3}(c_3)$. Therefore, we
have $\mathcal{R}^3_{\Phi}\subset\mathcal{S}^3_{\phi}$ where
\[
\mathcal{S}^3_{\phi}= \Biggl\{(f_1,f_2,f_3)\in
\mathcal{B}^3\dvtx\ell\Biggl(\sum_{l=1}^{3}f_{l} \Biggr)=1\mbox{ and,
for all } l\in\{1,2,3\}, \ell_{f_l}(c_{l})=y \Biggr\}.
\]
It follows that
\[
J(y)\geq\inf_{(f_1,f_2,f_3)\in\mathcal{S}^{3}_{\phi}}H \Biggl(\sum_{l=1}^{3}
f_{l}(x) \Biggr).
\]

\subsubsection*{Step 4: The simplified variational problem}

Recall that $y \in(c(B_1)/3, c(0)/3)$ is fixed in this part of
the proof. In this step, we prove that
%
%
\begin{equation}\label{eq:varpbI}
I(y):=\inf_{(f_1,f_2,f_3)\in\mathcal{S}^3_{\phi}}H \Biggl(\sum
_{l=1}^{3}f_{l}(x) \Biggr)
\end{equation}
is equal to $\Lambda^*(3y)$. Clearly, the set $\mathcal{S}^3_{\phi}$ is
convex. Therefore, if $\mathcal{S}^3_{\phi}$ is not empty, due to the
strict convexity of the relative entropy, the solution of the
variational problem (\ref{eq:varpbI}), say ${\mathbf{f}}^* =
(f_1^*,f_2^*,f_3^*) \in\mathcal{S}^3_{\phi}$, is unique, up to functions
which are null $\ell$-almost everywhere (a.e.). The variational problem
(\ref{eq:varpbI}) is an entropy maximization problem. We now compute
${\mathbf{f}}^*$ and check retrospectively that $\mathcal{S}^3_{\phi
}$ is not
empty. Consider the Lagrangian $\mathcal{L}$ defined by
\begin{eqnarray*}
&&
\mathcal{L}(f_1,f_2,f_3,\lambda_0,\lambda_1,\lambda_2,\lambda
_3)(x)\\
&&\qquad= \Biggl(\sum_{l=1}^3
f_l(x) \Biggr) \log\Biggl(\sum_{l=1}^3
f_l(x) \Biggr)+\lambda_0 \Biggl(\sum_{l=1}^3
f_l(x)- 1 \Biggr)\\
&&\qquad\quad{}+\sum_{l=1} ^3 \lambda_l\bigl(c_l(x)f_l(x)-y\bigr),
\end{eqnarray*}
where the $\lambda_i$'s $(i=0,\ldots,3)$ are the
Lagrange multipliers. For $l \in\{1,2,3\}$, define the Borel
sets
\[
A_l=\{x\in\mathbb{T}\dvtx f_l^*(x)>0\}.
\]
Since ${\mathbf{f}}^*$ is the solution of (\ref{eq:varpbI}), by
the Euler equations (see, e.g., Chapter 1 in \cite{buttazzo}) we
have, for $l\in\{1,2,3\}$,
\[
\biggl(\frac{\partial\mathcal{L}}{\partial
f_l} \biggr)\bigg|_{(f_1,f_2,f_3)={\mathbf{f}}^*}=0 \qquad\mbox{on }A_l.
\]
We deduce that, for all $x\in A_l$,
%
%
\begin{equation}\label{eq:diffl}
f_1^*(x)+f_2^*(x)+f_3^*(x)=\rme^{-1-\lambda_0-\lambda_l
c_l(x)}.
\end{equation}
Define the functions $g_{1}(x):=f^{*}_{2}(jx)$,
$g_{2}(x):=f^{*}_{3}(jx)$ and $g_{3}(x):=f^{*}_{1}(jx)$. By a
change of variable, it is straightforward to check that
$(g_1,g_2,g_3)\in\mathcal{S}^{3}_{\phi}$ and
\[
\int_{\mathbb{T}} \Biggl(\sum_{l=1}^{3}g_{l}(x) \Biggr)
\log\Biggl(\sum_{l=1}^{3}g_{l}(x) \Biggr) \,dx=\int_{\mathbb{T}}
\Biggl(\sum_{l=1}^{3}f^{*}_{l}(x) \Biggr)\log\Biggl(\sum_{l=1}^3
f^{*}_{l}(x) \Biggr) \,dx.
\]
The uniqueness of the solution implies that a.e.
\[
f^{*}_{2}(jx)=f^{*}_{1}(x),\qquad f^{*}_{3}(jx)=f^{*}_{2}(x)\quad\mbox{and}
\quad f^{*}_{1}(jx)=f^{*}_{3}(x).
\]
In particular, up to a null measure set, $A_{l}=j^{l-1}A_1$.
Moreover, on $A_1$, the equality, a.e. $\sum_{l=1}^{3}g_{l}(x)=
\sum_{l=1}^{3}f^{*}_{l}(x)$ applied to (\ref{eq:diffl}) gives,
a.e. on $A_1$, $\exp(-1-\lambda_{0}-\lambda_{2}c_2(jx))=\exp(-1-
\lambda_{0}-\lambda_{1}c_1(x))$ (indeed $x\in A_{1}$ implies $jx
\in A_2$). We deduce that $\lambda_{2}=\lambda_{1}$. The same
argument on $A_3$ carries over by symmetry, so finally $\lambda_1
=\lambda_2=\lambda_3$. We now
use the following lemma that will be proved at the end of the
step.
\begin{Lemma} \label{le:B1T1}
Under the foregoing assumptions and notation, up to a Borel set of
null Lebesgue measure it holds $A_1\subset\mathbb{T}_1$.
\end{Lemma}

By Lemma \ref{le:B1T1} and the a.e. equality $A_{l}=j^{l-1}A_1$,
we deduce that $A_{1}\in\mathcal{T}$, up to a Borel set of null Lebesgue
measure. So, by (\ref{eq:diffl}) and the equality
$\lambda_1=\lambda_2=\lambda_3$, it follows that
\[
f_1^*(x)=\rme^{-1-\lambda_{0}-\lambda_{1}c(x)}\mathbh{1}(x\in
A_1)\qquad \mbox{a.e.}
\]
and $f^*_2(x)=f^*_1(j ^2 x)$, $f^*_3(x)=f^*_1(j
x)$. Note that the constraints
\[
\ell\Biggl(\sum_{l=1}^{3}f_{l}^{*} \Biggr)=1 \quad\mbox{and}\quad
\ell_{f^*_1}(c_{1})=y
\]
read, respectively,
\[
\int_{A_1}\rme^{-1-\lambda_{0}-\lambda_{1}c(x)} \,dx
=1/3 \quad\mbox{and}\quad \int_{A_1}c(x)\rme^{-1-
\lambda_{0}-\lambda_{1}c(x)} \,dx=y.
\]
This implies that the Lagrange multipliers $\lambda_0$ and
$\lambda_1$ are solutions of the equations in
(\ref{eq:gammaeq0}) with $B=A_1$. Moreover
\begin{eqnarray*}
\int_{\mathbb{T}} \Biggl(\sum_{l=1}^{3}f^{*}_{l}(x) \Biggr)
\log\Biggl(\sum_{l=1}^{3}f^{*}_{l}(x) \Biggr) \,dx&=&3
\int_{A_1} \bigl(-1-\lambda_{0}-\lambda_{1}c(x) \bigr)
\rme^{-1-\lambda_0-\lambda_{1}c(x)} \,dx\\
&=&-(1+\lambda_0)-3y\lambda_1.
\end{eqnarray*}
Therefore (see the beginning of step 2)
\[
I(y)=\int_{\mathbb{T}} \Biggl(\sum_{l=1}^{3}f^{*}_{l}(x) \Biggr)
\log\Biggl(\sum_{l=1}^{3}f^{*}_{l}(x) \Biggr) \,dx=
\nu(A_1).
\]
Since $A_{1}\in\mathcal{T}$ we deduce that
\[
I(y)\geq\inf\{\nu(B)\dvtx B\in\mathcal{T}\}.
\]
For the reverse inequality, take $B\in\mathcal{T}$ such that
$\nu(B)=\sup_{(\eta_0,\eta_1)\in\mathbb{R}^2}F(B,\eta_0,\eta
_1)$ is
finite. Since the function $(\eta_0,\eta_1)\mapsto F(B,
\eta_0,\eta_1)$ is finite and strictly concave, it admits a unique
point of maximum. Arguing exactly as at the beginning of step~2,
we have that the point of maximum is $(\gamma_0(B),\gamma_1(B))$,
whose components are solutions of equations in
(\ref{eq:gammaeq0}), and
\[
\nu(B)=-\bigl(1+\gamma_0(B)\bigr)-3y\gamma_1(B).
\]
For $l\in\{1,2,3\}$, define the functions on $\mathbb{T}$
\[
g_{l,B}\dvtx x\mapsto\rme^{-1-\gamma_0(B)-\gamma_1(B)c_{l}(x)}
\mathbh{1}(x\in j^{l-1}B).
\]
Since $\gamma_0(B)$ and $\gamma_1(B)$ solve the equations in
(\ref{eq:gammaeq0}), it follows easily that
$(g_{1,B}, g_{2,B},g_{3,B})\in\mathcal{S}^3 _\phi$. Therefore
\begin{eqnarray*}
\nu(B)&=&\int_{\mathbb{T}} \Biggl(\sum_{l=1}^{3}g_{l,B}(x) \Biggr)
\log\Biggl(\sum_{l=1}^{3}g_{l,B}(x) \Biggr) \,dx\\
&\geq&\inf_{(f_1,f_2,f_3)\in\mathcal{S}^3_{\phi}}H \Biggl(\sum_{l=1}^{3}
f_{l}(x) \Biggr).
\end{eqnarray*}
Thus
\[
I(y)=\nu(A_1)=\inf\{\nu(B)\dvtx B\in\mathcal{T}\}.
\]
Since $A_1\in\mathcal{T}$, by Lemma \ref{le:minmax} we get that
$\ell(A_1)=1/3$. So, by Lemma \ref{le:B1T1}, we deduce that
$A_{1}=\mathbb{T}_1$ up to a Borel set of null Lebesgue measure. Then by
(\ref{eq:nuT1}) we conclude
\[
I(y)=\Lambda^{*}(3y).
\]
\begin{pf*}{Proof of Lemma \ref{le:B1T1}}
The argument is by contradiction. Define the Borel set
\[
C:=(A_{1}\cap\mathbb{T}_{1}^{c})\cup(jA_{1}\cap\mathbb
{T}_{2}^{c})\cup(j^2A_{1}
\cap\mathbb{T}_{3}^{c})
\]
and assume that $\ell(A_{1}\cap\mathbb{T}_{1}^{c})>0$. For
$l\in\{1,2,3\}$, define $\widetilde{A}_{l}=(A_{l}\setminus
C)\cup(C\cap\mathbb{T}_l)$ and $\widetilde{g}_{l}(x)=(f^{*}_1(x)+
f^{*}_2(x)+f^{*}_3(x))\mathbh{1}(x\in\widetilde{A}_l)$. Since
$A_{l}=j^{l-1}A_1$ up to a Borel set of null Lebesgue measure,
then $j^{l-1}C=C$ and $\widetilde{A}_{l}=j^{l-1}\widetilde{A}_1$
up to a Borel set of null Lebesgue measure. So by
(\ref{eq:diffl}) it follows that
$\ell_{\widetilde{g}_1}(c_{1})=\ell_{\widetilde{g}_{2}}(c_2)=
\ell_{\widetilde{g}_3}(c_{3})$, and therefore
%
%
\begin{equation}\label{eq:eqwith3}
3\int_{\mathbb{T}}c_{l}(x)\widetilde{g}_{l}(x) \,dx=\int_{\mathbb{T}}
\Biggl(\sum_{l=1}^{3}\mathbh{1}(x\in\widetilde{A}_{l})c_{l}(x) \Biggr)
\Biggl(\sum_{l=1}^{3}f^{*}_{l}(x) \Biggr) \,dx.
\end{equation}
Now, note that $\widetilde{A}_{l}\subseteq\mathbb{T}_l$ and,
up to a Borel set of null Lebesgue measure,
%
%
\begin{equation}\label{eq:eqaatilde}
\widetilde{A}_{1}\cup\widetilde{A}_{2}\cup\widetilde{A}_{3}=A_{1}
\cup A_{2}\cup A_{3}.
\end{equation}
So by assumption (\ref{ass3}), a.e.
\[
\mathbh{1}(x\in\widetilde{A}_{l})c_{l}(x)\leq\sum_{m=1}^{3}\mathbh
{1}(x\in
A_m)c_m(x),
\]
and the inequality is strict if $x$ is in
$C\cap\accentset{\circ}{\mathbb{T}}_l$. Indeed if $x\in
C\cap\accentset{\circ}{\mathbb{T}}_l$, then a.e. $x\in A_{m}$ for some
$m\ne
l$, and so $c_{l}(x)<c_{m}(x)$ by (\ref{ass3}). Therefore, since
$\ell(A_1\cap\mathbb{T}_1^c)>0$ then $\ell(C\cap\accentset{\circ
}{\mathbb{T}}_l)>0$
and, using (\ref{eq:eqwith3}), we get
\[
\int_{\mathbb{T}}c_{l}(x)\widetilde{g}_{l}(x) \,dx<\int_{\mathbb{T}
}c_{l}(x)f^{*}_{l}
(x) \,dx=y.
\]
For $p\in[0,1]$, define the functions
\[
\widetilde{g}_{l,p}(x)=(1-p) \widetilde{g}_{l}(x)+p \mathbh{1}\bigl(x\in
\mathbb{T}
_{\sigma(l)}\bigr),
\]
where $\sigma= (1\ 2\ 3)$ is the cyclic permutation.
By assumption (\ref{eq:concave1}) it follows that
\[
\int_{\mathbb{T}}c_{l}(x)\widetilde{g}_{l,1}(x) \,dx>c(0)/3>y.
\]
We have already checked that
$\ell_{\widetilde{g}_{l,0}}(c_l)< y$, thus, by the mean value
theorem, there exists $\overline{p}\in(0,1)$ such that
$(\widetilde{g}_{1,\overline{p}},\widetilde{g}_{2,\overline{p}},\widetilde
{g}_{3,\overline{p}})\in\mathcal{S}^3_{\phi}$. The
convexity of the relative entropy gives
\begin{eqnarray*}
H(\widetilde{g}_{1,\overline{p}}+\widetilde{g}_{2,\overline{p}}+\widetilde
{g}_{1,\overline{p}} | \ell)&\leq&
\overline{p}H(\widetilde{g}_1+\widetilde{g}_2+\widetilde{g}_3 | \ell)+
(1-\overline{p})H(\ell| \ell)\\
&=&\overline{p}H(f^{*}_{1}+f^{*}_{2}+{f}^{*}_3
|
\ell),
\end{eqnarray*}
where the latter equality follows by
(\ref{eq:eqaatilde}) and the definition of
$\widetilde{g}_l$. This leads to a contradiction since ${\mathbf{f}}=
(f^*_1,f^*_2,f^*_3)$ minimizes the relative entropy on
$\mathcal{S}^3_{\phi}$.
\end{pf*}

\subsubsection*{Step 5: End of the proof}

It remains to check that ${\mathbf{f}}^{*}=(f^*_1,f^*_2,f^*_3)\in
\mathcal{R}^{3}_{\Phi}$. For this we need to prove that
$\Phi(\ell_{f^*_1+f^*_2+f^*_3})=\phi(\ell_{f^*_1},\ell
_{f^*_2},\ell_{f^*_3})
=y$. Since ${\mathbf{f}}^{*}\in\mathcal{S}^{3}_{\phi}$ then
$\ell_{f_1^*}(c_1)=\ell_{f_2^*}(c_2)=\ell_{f_3^*}(c_3)=y$;
moreover, by the properties of the functions $f_l^*$ it holds
$\ell_{f_l^*}(c_l)=\int_{\mathbb{T}_l}c_l(x)f_l(x) \,dx$. So the
claim follows if we check that
\[
\Phi(\ell_{f^*_1+f^*_2+f^*_3})\geq\int_{\mathbb{T}_1}c_1(x)f_1(x)
\,dx.
\]
By Lemma \ref{le:tilderho}(i) we have that there\vspace*{1pt}
exists $(g_1,g_2,g_3)\in\mathcal{B}^3$ such that
$\ell_{f^*_1+f^*_2+f^*_3}=\ell_{g_1}+\ell_{g_2}+\ell_{g_3}$,
$\Phi(\ell_{f^*_1+f^*_2+f^*_3})=\phi(\ell_{g_1},\ell_{g_2},\ell_{g_3})$
and $\ell_{g_1}(c_1)=\ell_{g_2}(c_2)=\ell_{g_3}(c_3)$. In
particular,
%
%
\begin{eqnarray}\label{eq:cond3}\qquad
3\Phi(\ell_{f^*_1+f^*_2+f^*_3}
)&=&\sum_{l=1}^{3}\int_{\mathbb{T}}c_{l}(x)g_{l}(x) \,dx
=\sum_{m=1}^{3}\int_{\mathbb{T}_m}\sum_{l=1}^{3}c_{l}(x)g_{l}(x) \,
dx\nonumber\\
&\geq&\sum_{m=1}^{3}\int_{\mathbb{T}_m}c_{m}(x)\sum_{l=1}^{3}g_{l}(x)
\,dx\nonumber\\[-8pt]\\[-8pt]
&\geq&\sum_{m=1}^{3}\int_{\mathbb{T}_m}c_{m}(x)f_{m}^{*}(x) \,
dx\nonumber\\
&=&3\int_{\mathbb{T}_1}c_{1}(x)f_{1}^{*}(x) \,dx,\nonumber
\end{eqnarray}
where in (\ref{eq:cond3}) we used assumption
(\ref{ass3}). This concludes the proof of Theorem~\ref{th:LDP}(i).

\subsection{\texorpdfstring{Proof of Theorem
\protect\ref{th:LDP}\textup{(ii)}}{Proof of Theorem 1.3(ii)}}

Some ideas in the following proof of Theorem~\ref{th:LDP}(ii)
are similar to those one in the proof of Theorem
\ref{th:LDP}(i). Therefore, we shall omit some details. We
divide the proof of Theorem \ref{th:LDP}(ii) in 3 steps.

\subsubsection*{Step 1: Case $y\protect\notin(c(B_1)/3,c(0))$}

As noticed in step 1 of the proof of Theorem~\ref{th:LDP}(i),
for any measure $\beta\in\mathcal{M}_b(\mathbb{T})$, $\beta
(c_l)\geq
c(B_1)\beta(\mathbb{T})$, and the equality holds only if
$\beta=\delta_{B_l}$. We deduce that, for all
$\alpha\in\mathcal{M}_1(\mathbb{T})$, $3\Psi(\alpha)>c(B_1)$.
Therefore, by
Theorem \ref{th:LDPc}(ii), $\overline{J}(y)=+\infty$ if $y\leq
c(B_1)/3$. Now, note that, for $\alpha\in\mathcal{M}_1(\mathbb{T})$
it holds that
\[
\Psi(\alpha)=\max_{1\leq l\leq3} \biggl(\int_{\mathbb
{T}_l}c_{l}(x)\alpha
(dx) \biggr)<c(0)\max_{1\leq l\leq3}\alpha(\mathbb{T}_l)\leq
c(0),
\]
where the strict inequality follows by assumption
(\ref{eq:concave}) and $\alpha\ll\ell$. Therefore, using again
Theorem \ref{th:LDPc}(ii), we easily deduce that
$\overline{J}(y)=+\infty$ if $y\geq c(0)$.

\subsubsection*{Step 2: The set function $\mu$}

For the remainder of the proof we fix $y\in(c(B_1)/3,c(0))$, and
we shall often omit the dependence on $y$ of the quantities under
consideration. In the following we argue as in step 2 of the proof
of Theorem \ref{th:LDP}(i). Let $B\subset\mathbb{T}$ be a Borel set
with positive Lebesgue measure and define the function of
$(\eta_0,\eta_1)\in\mathbb{R}^2$
\[
q(B,\eta_0,\eta_1)=2\rme^{-1-\eta_0}\ell(B\cap\mathbb{T}_2)+\int
_{B\cap
\mathbb{T}_1}\rme^{-1-\eta_{0}-\eta_{1}c(x)}
\,dx.
\]
Clearly, $q(B,\cdot)$ is strictly convex on
$\mathbb{R}^2$. Define the strictly concave function
\[
G(B,\eta_0,\eta_1)=-\eta_0- y\eta_1-q(B,\eta_0,\eta_1)
\]
and the set function
\[
\mu(B)=\sup_{(\eta_0,\eta_1)\in\mathbb{R}^2}G(B,\eta_0,\eta_1).
\]
If there exist $\overline{\gamma}_0=\overline{\gamma}_{0}(B)$ and
$\overline{\gamma}_1=\overline{\gamma}_{1}(B)$
such that
%
%
\begin{eqnarray}\label{eq:gammaeq1}
\int_{B\cap\mathbb{T}_1} \rme^{-\overline{\gamma}_{1}c(x)}\,
dx+2\ell(B
\cap\mathbb{T}_2)& = & \rme^{1+\overline{\gamma}_{0}}
\quad\mbox{and}\nonumber\\[-8pt]\\[-8pt]
\int_{B
\cap\mathbb{T}_1} c(x)\rme^{-\overline{\gamma}_{1}
c(x)}\,dx &=& y\rme^{1+\overline{\gamma}_{0}}, \nonumber
\end{eqnarray}
then we have
\[
\mu(B)=-\bigl(1+\overline{\gamma}_{0}(B)\bigr)-y\overline{\gamma}_{1}(B).
\]
In particular, by Proposition \ref{prop:Laplace}(ii),
setting $\overline{\gamma}_{1}(\mathbb{T})=-\eta_y$ and $\overline
{\gamma}_{0}(\mathbb{T})=\overline{\Lambda}(\eta_y)-1$
one has
%
%
\begin{equation}\label{eq:oL*y}\quad
\overline{\Lambda}{}^*(y)=\mu(\mathbb{T})=-\bigl(1+\overline{\gamma
}_{0}(\mathbb{T})\bigr)-y\overline{\gamma}_{1}(\mathbb{T})\qquad
\mbox{if }\gamma<y<c(0),
\end{equation}
and $\overline{\gamma}_{0}(\mathbb{T})$ and $\overline{\gamma
}_{1}(\mathbb{T})$ are the unique
solutions of the equations in (\ref{eq:gammaeq1}) with
$B=\mathbb{T}$. Recall also that in step 2 of the proof of Theorem
\ref{th:LDP}(i) we showed
\[
\Lambda^*(3y)=-\bigl(1+\gamma_{0}(\mathbb{T}_1)\bigr)-3y\gamma_{1}(\mathbb{T}_1)\qquad
\mbox{if }c(B_1)/3<y\leq\gamma,
\]
where $\gamma_{0}(\mathbb{T}_1)$ and $\gamma_{1}(\mathbb{T})$ are the
unique solutions of the equations in (\ref{eq:gammaeq0}) with
$B=\mathbb{T}_1$. Note that, for Borel sets $A$ and $B$ such that
$A\subseteq B\subseteq\mathbb{T}$, we have, for all $\eta_0,\eta
_1\in\mathbb{R}$,
$G(A,\eta_0,\eta_1)\geq G(B,\eta_0,\eta_1)$. This proves that the
set function $\mu$ is nonincreasing (for the set inclusion). An
easy consequence is the following lemma:
\begin{Lemma}\label{le:minmax1}
Under the foregoing assumptions and notation, it holds that
\[
\inf\{\mu(B)\dvtx B\subseteq\mathbb{T}\}=\overline{\Lambda}{}^*(y)
\qquad\mbox{if }\gamma<y<c(0).
\]
\end{Lemma}

\subsubsection*{Step 3: The related variational problem}

As above we fix $y\in(c(B_1)/3,c(0))$; as in the proof of Theorem
\ref{th:LDP}(i) we denote by $\mathcal{B}$ the set of Borel
functions defined on $\mathbb{T}$ with values in $[0,\infty)$. By Theorem
\ref{th:LDPc}(ii), we have
\[
\overline{J}(y)=\inf_{f\in\mathcal{U}}H(f),
\]
where
\[
\mathcal{U}= \biggl\{f\in\mathcal{B}\dvtx\ell(f)=1 \mbox{ and } \max
_{1\leq l
\leq3} \biggl(\int_{\mathbb{T}_l}c_l(x)f(x)\,dx \biggr)=y \biggr\}.
\]
Note that $f\in\mathcal{U}$ if and only if the functions $x\mapsto f(jx)$
and $x\mapsto f(j^{2}x)$ are also in $\mathcal{U}$ and so
%
%
\begin{equation}\label{eq:infoJ}
\overline{J}(y)=\inf_{f\in\mathcal{V}}H(f),
\end{equation}
where
\[
\mathcal{V}= \{f\in\mathcal{B}\dvtx\ell(f)=1, {\ell
_{f}}_{|\mathbb{T}_1}(c_1)=y,
{\ell_{f}}_{|\mathbb{T}_2}(c_2)\leq y, {\ell_{f}}_{|\mathbb
{T}_3}(c_3)\leq y
\}.
\]
The optimization problem (\ref{eq:infoJ}) is a
minimization of a convex function on a convex set defined by
linear constraints. Thus it can be solved explicitly. Therefore,
if $\mathcal{V}$ is not empty, since the relative entropy is strictly
convex, the solution of the variational problem
(\ref{eq:infoJ}), say $f^{*}\in\mathcal{V}$, is unique, up to functions
which are null $\ell$-almost everywhere. We will compute $f^*$ and
show that $\mathcal{V}$ is not empty at the same time. So assume that
$\mathcal{V}$ is not empty and define the function
\[
g(x)=f^{*}(x)\mathbh{1}_{\mathbb{T}_1}(x)+f^{*}(jx)\mathbh
{1}_{\mathbb{T}_2}(x)+f^{*}(j^{2}
x)\mathbh{1}_{\mathbb{T}_3}(x).
\]
It is easily checked that $g\in\mathcal{V}$ and $H(g)=H(f)$. The
uniqueness of $f^*$ implies that
%
%
\begin{equation}\label{eq:T2idT3}
\mbox{for almost all }x\in\mathbb{T}_2\qquad f^*(jx)=f^*(x).
\end{equation}
Therefore, up to modifying $f^*$ on a set of null
measure, $f^*\in\mathcal{V}'$ where
\[
\mathcal{V}'= \{f\in\mathcal{B}\dvtx\ell(f)=1, {\ell
_{f}}_{|\mathbb{T}_1}(c_1)=y,
{\ell_{f}}_{|\mathbb{T}_2}(c_2)\leq y \}
\]
and the variational problem reduces to
$\overline{J}(y)=\inf_{f\in\mathcal{V}'}H(f)$. Consider the
Lagrangian $\mathcal{L}$
defined by
\begin{eqnarray*}
\mathcal{L}(f,\lambda_0,\lambda_1,\lambda_2)(x)&=&f(x)\log
f(x)+\lambda_{0}\bigl(f(x)-1\bigr)+\lambda_1\bigl(c_1(x)f(x)\mathbh{1}_{\mathbb{T}_1}(x)-y\bigr)\\
&&{}+\lambda_{2}\bigl(c_2(x)f(x)\mathbh{1}_{\mathbb{T}_2}(x)-y\bigr)
\end{eqnarray*}
with
\[
\lambda_{2} \biggl(\int_{\mathbb{T}_2}c_{2}(x)f^{*}(x) \,dx-y \biggr)
=0.
\]
The two cases $\lambda_{2}=0$ (i.e., $f^*$ is not
constrained on $\mathbb{T}_2$) and $\lambda_{2}\ne0$ (i.e., $f^*$ is
constrained on $\mathbb{T}_2$) are treated separately. For each case, we
solve the variational problem. The optimal function is denoted by
$f_u$ for $\lambda_{2}=0$ and by $f_c$ for $\lambda_{2}\ne0$, so
that $f^{*}=\arg\min(H(f_u),H(f_c))$. Assume first that
$\lambda_{2}=0$ so that $f^* = f_u$ and define the Borel set
\[
A_{u}=\{x\in\mathbb{T}\dvtx f_u(x)>0\}.
\]
By the Euler equations (see, e.g., Chapter 1 in
\cite{buttazzo}) we get, for all $x\in\mathbb{T}$,
%
%
\begin{equation}\label{eq:fu}
f_u(x)=\mathbh{1}_{\mathbb{T}_1\cap A_u}(x)\rme^{-1-\lambda
_0-\lambda_1
c_1(x)}+\mathbh{1}_{(\mathbb{T}_2\cup\mathbb{T}_3)\cap A_u}(x)
\rme^{-1-\lambda_0}.
\end{equation}
By (\ref{eq:T2idT3}) we have
$\ell(A_u\cap\mathbb{T}_2)=\ell(A_u\cap\mathbb{T}_3)$, and so the
constraints
$\ell(f_u)=1$ and ${\ell_{f_u}}_{|\mathbb{T}_1}(c_1)=y$ read,
respectively,
\[
\int_{A_{u}\cap\mathbb{T}_1}\rme^{-\lambda_{1}c(x)} \,dx
+2\ell(A_{u}\cap\mathbb{T}_{2})=\rme^{1+\lambda_0}
\]
and
\[
\int_{A_{u}\cap\mathbb{T}_1} c(x) \rme^{-\lambda_1
c(x)} \,dx=y\rme^{1+\lambda_0}.
\]
%
With the notation of step 2, this implies that
$\lambda_0=\overline{\gamma}_0(A_u)$ and $\lambda_1=\overline
{\gamma}_1(A_u)$ are the solution
of the equations in (\ref{eq:gammaeq1}) with $B=A_u$. In
particular,
\[
\mu(A_u)=-\bigl(1+\overline{\gamma}_0(A_u)\bigr)-y\overline{\gamma}_1(A_u)=H(f_u),
\]
where the latter equality follows from the computation
of the entropy using (\ref{eq:fu}). By Lemma \ref{le:minmax1} we
deduce that
\[
H(f_u)\geq\overline{\Lambda}{}^*(y) \qquad\mbox{if }\gamma<y<c(0).
\]
By (\ref{eq:oL*y}) we have $H(h)=\overline{\Lambda}{}^*(y)$, where
\[
h(x)=\mathbh{1}_{\mathbb{T}_1}(x)\rme^{-1-\overline{\gamma}_{0}
(\mathbb{T}) - \overline{\gamma}_{1} (\mathbb{T})
c(x)}+\mathbh{1}_{\mathbb{T}_2\cup\mathbb{T}_3}(x)\rme
^{-1-\overline{\gamma}_{0} (\mathbb{T})},
\]
and $\overline{\gamma}_{0}(\mathbb{T})$, $\overline{\gamma
}_{1}(\mathbb{T})$ are the unique
solutions of the equations in (\ref{eq:gammaeq1}) with $B=\mathbb{T}$.
Now we prove that $h\in\mathcal{V}$, for $\gamma<y<c(0)$, so that
%
%
\begin{equation}\label{eq:2=0}
H(f_u)=\overline{\Lambda}{}^*(y) \qquad\mbox{if }\gamma<y<c(0).
\end{equation}
Recall that $-\overline{\gamma}_1(\mathbb{T})$ is the unique
solution of
\[
\frac{\int_{\mathbb{T}_1}c(x)\rme^{\theta c(x)} \,dx}
{\int_{\mathbb{T}_1}\rme^{\theta c(x)} \,dx+2/3}=y.
\]
The function
\[
\theta\mapsto\frac{\int_{\mathbb{T}_1}c(x)\rme^{\theta
c(x)} \,dx}{\int_{\mathbb{T}_1}\rme^{\theta
c(x)} \,dx+ 2/3}
\]
is strictly increasing (as can be checked by a straightforward
computation) and, for $\theta=0$, it is equal to
$\gamma$. Therefore, since $y>\gamma$, we have $-\overline{\gamma
}_1(\mathbb{T})>0$. It
implies that
\[
\int_{\mathbb{T}_1}c(x)\rme^{-1-\overline{\gamma}_0(\mathbb
{T})-\overline{\gamma}_1(\mathbb{T})c(x)}
\,dx
=y>\int_{\mathbb{T}_1}c(x)\rme^{-1-\overline{\gamma}_0(\mathbb
{T})} \,dx=\gamma
\rme^{-1-\overline{\gamma}_0(\mathbb{T})}.
\]
In particular, $h\in\mathcal{V}$. Now we deal with the case
$\lambda_{2}\neq0$. We have
\[
{\ell_{f_c}}_{|\mathbb{T}_1}(c_{1})={\ell_{f_c}}_{|\mathbb
{T}_2}(c_{2})={\ell
_{f_c}}_{|\mathbb{T}_3}(c_{3})
= y.
\]
In particular, if we set $f_{c,l}(x)=\mathbh{1}(x\in\mathbb{T}_{l})
f_{c}(x)$, we get $(f_{c,1},f_{c,2},f_{c,3})\in\mathcal{S}^{3}_\phi
$. By
step 4 of the proof of Theorem \ref{th:LDP}(i), it implies that
\[
H(f_c)\geq\inf_{(f_{1},f_{2},f_{3})\in\mathcal{S}^3_\phi}H(f_{1}+f_{2}+
f_{3})=\Lambda^{*}(3y)=H(f_{1}^{*}+f_{2}^{*}+f_{3}^{*}),
\]
where ${\mathbf{f}}^{*}=(f_1^*,f_2^*,f_3^*)$ was defined above. Since
$f_1^* + f_2 ^* + f_3 ^* \in\mathcal{V}$, we deduce directly that a.e.
$f_c = f_1^* + f_2 ^* + f_3 ^*$ and
%
%
\begin{equation}\label{eq:2ne0}
H(f_c)=\Lambda^*(3y).
\end{equation}
It remains to find out for which values of $y$ the
Lagrange multiplier $\lambda_2$ is equal to zero. First of all
note that if $y=\gamma$, then the function identically equal to $1$
is in $\mathcal{V}$. We deduce that $f^*\equiv1$ and so $\lambda_2=0$
(since the optimal solution is not constrained on $\mathbb{T}_2$) and
$\overline{J}(\gamma)=0=\Lambda^*(3\gamma)$. Now assume $\gamma
<y<c(0)$. By
Proposition \ref{prop:Laplace}(iii), we deduce
$\overline{\Lambda}{}^*(y)<\Lambda^*(3y)$. It follows by (\ref
{eq:2=0}) and
(\ref{eq:2ne0}) that $H(f_u)<H(f_c)$. Recall that $f^{*}=
\arg\min(H(f_u),H(f_c))$, thus $\lambda_{2}=0$ and
$\overline{J}(y)=\overline{\Lambda}{}^*(y)$. It remains to deal with
the case
$c(B_1)/3<y<\gamma$. The following lemma holds:
\begin{Lemma} \label{le:48}
Under the foregoing assumptions and notation, if
$c(B_1)/3<y<\gamma$, then $\overline{J}(y)\geq J(y)$.
\end{Lemma}

Then, by Theorem \ref{th:LDP}(i) and (\ref{eq:2ne0}) we get
\[
\Lambda^{*}(3y) = J(y)
\leq\overline{J}(y)=\min(H(f_u),H(f_c))\leq\Lambda^*(3y).
\]
This completes the proof.
\begin{pf*}{Proof of Lemma \ref{le:48}} Choose
$y<z<\gamma$. By construction $P(\overline{\rho}_n\leq nz)\leq
P(\rho_{n}\leq nz)$. Taking the logarithm, applying Theorem
\ref{th:LDPc} and recalling that $\overline{J}(y)=J(y)=+\infty$ for
$y\leq
c(B_1)/3$ we have
\begin{eqnarray*}
-\inf_{t\in(c(B_1)/3,z)}\overline{J}(t)&\leq&\liminf_{n\to\infty
}\frac
{1}{n}\log
P(\overline{\rho}_{n}\leq nz) \\
&\leq& \limsup_{n\to\infty}\frac{1}{n}
\log P(\rho_{n}\leq nz)\\
&\leq&-\inf_{t\in(c(B_1)/3,z]}J(t).
\end{eqnarray*}
Therefore
\[
\overline{J}(y)\geq\inf_{t\in(c(B_1)/3,z)}\overline{J}(t)\geq\inf
_{t\in
(c(B_1)/3,z]}J(t)=J(z),
\]
where the latter equality follows since $J(y)=\Lambda^*(3y)$ is
decreasing on $(c(B_1)/\break3,\gamma)$. Recalling that $J(y)
=\Lambda^*(3y)$ is also continuous on $(c(B_1)/3,\gamma)$, the
claim follows letting $z$ tend to $y$.
\end{pf*}

\section{Model extension}\label{sec:extension}

\subsection{The analog one-dimensional model}

The analog one-dimensional model is obtained as follows. There are
$n$ objects on $(0,1)$, say $\{1,\ldots,n\}$, and two bins located
at 0 and 1, respectively. The location of the $k$th object is
given by a r.v. $X_k$ and it is assumed that the r.v.'s
$\{X_k\}_{1\leq k \leq n}$ are i.i.d. and uniformly distributed on
$[0,1]$. The cost to allocate an object at $x\in[0,1]$ to the bin
at $0$, respectively, at $1$, is $c(x)$, respectively, $c(1-x)$. The
asymptotic analysis of allocations which realize the optimal and
the suboptimal load can be carried on using the ideas and the
techniques developed in this paper. Due to the simpler geometry of
the one-dimensional model, many technical difficulties met in the
two-dimensional case disappear, and with the proper assumptions on
the cost function, it is possible to state and prove the analog of
Theorems \ref{th:lln}, \ref{th:tcl} and \ref{th:LDP}.

\subsection{Random cost function}
An interesting and natural extension of the model takes into
account random cost functions. Let $\mathcal{Z}$ be a Polish space and
$\mathbf{Z}_{k}=(Z_{k}^1,Z_{k}^{2},Z_{k}^{3})$ ($k=1,\ldots,n$) a
r.v. taking values on $\mathcal{Z}^3$. Assume that: the sequences
$\{X_k\}_{1\leq k\leq n}$ and $\{\mathbf{Z}_k\}_{1\leq k\leq n}$ are
independent; the r.v.'s $\{\mathbf{Z}_k\}_{1\leq k\leq n}$ are
i.i.d. with common distribution $Q$; the r.v.'s $Z^1_1$, $Z^2_1$
and $Z^3_1$ are i.i.d. Let $c\dvtx\mathbb{T}\times\mathcal{Z}^3\to
[0,\infty)$ be a
measurable function. We consider an extension of the basic model
where the cost to allocate the $k$th object to the bin at $B_l$
($l=1,2,3$) is equal to $c_l(X_k,\mathbf{Z}_{k})$. Here, for
$\mathbf{z}=(z^1,z^2,z^3)$, the cost functions are defined in such a
way that they preserve the spatial symmetry
$c_1(x,\mathbf{z})=c(x,\mathbf{z})$, $c_{2}(x,\mathbf{z})=c
(j^{2}x,(z^2,z^3,z^1))$ and
$c_3(x,\mathbf{z})=c(jx,(z^3,z^1,z^2))$. The load associated to an
allocation matrix $A\in\mathcal{A}_n$ is
\[
\rho_{n}(A)=\max_{1\leq l\leq3} \Biggl(\sum_{k=1}^{n}a_{kl}c_l
(X_k,\mathbf{Z}_k) \Biggr).
\]
In a wireless communication scenario we have $\mathcal{Z}=\mathbb
{R}_+$, and the
typical cost function is of the form
\[
c(x,\mathbf{z})=\frac{a+\min\{b,z^{2}|x-B_2|^{-\alpha}\}+\min\{b,
z^{3}|x-B_3|^{-\alpha}\}}{\min\{b,z^{1}|x-B_1|^{-\alpha}\}},
\]
where $a>0$, $\alpha\geq2$ and
$b>(\lambda\sqrt{3}/2)^{-\alpha}$. The additional randomness in
the cost function models the fading along the channel (see, e.g.,
\cite{tse}). The suboptimal allocation $\overline{A}=(
\overline{a}_{k,l})_{1\leq k\leq n,1\leq l\leq3}$ is obtained by
allocating each point to its less costly bin. To be more precise,
assume that $\ell\otimes Q$-a.s., for any $l\neq m$, $c_{l}(x,\mathbf{z})
\neq c_{m}(x,\mathbf{z})$. Then,
setting
\[
\overline{a}_{k,l}=\mathbh{1}\Bigl(c_{l}(X_k,\mathbf{Z}_k)<\min_{m \neq l}c_{m}(
X_k,\mathbf{Z}_k)\Bigr),
\]
the suboptimal allocation matrix is a.s. well defined. Consider
the suboptimal load $\overline\rho_{n}=\rho_{n}(\overline{A})$ and
the optimal load $\rho_{n}=\min_{A\in\mathcal{A}_n}\rho_{n}(A)$. Exactly
as in the proof of Theorem \ref{th:lln}, one can prove that, a.s.
\[
\lim_{n\to\infty}\frac{\rho_{n}}{n}=\lim_{n\to\infty}
\frac{\overline\rho_n}{n}=\int_{\mathbb{T}\times\mathcal
{Z}^3}\mathbh{1}
\Bigl(c_{l}(x,\mathbf{z})
<\min_{m \neq l}c_{m}(x,\mathbf{z})\Bigr) \,dxQ(d\mathbf{z}).
\]
Deriving analogs of Theorem \ref{th:tcl} and Theorem \ref{th:LDP} is an
interesting issue. For the central limit theorem, an analog
of the suboptimal allocation matrix $\hat{A}$ in Proposition
\ref{prop:subtcl} should be defined. For the large deviation
principles, the contraction principle can be applied as well, but
it might be more difficult to solve the associated variational
problems.

\subsection{Asymmetric models}

Most techniques of the present paper collapse when the symmetry of
the model fails, for example, the region is not an equilateral triangle,
the locations are not uniformly distributed on the triangle, the
cost of an allocation is not properly balanced among the bins. For
a result on the law of large numbers in the case of an asymmetric
model, we refer the reader to Bordenave \cite{CBIEEE06}.

\begin{appendix}\label{sec:app}
\section*{Appendix}

\subsection{\texorpdfstring{Proof of Lemma \protect\ref{le:tvd}}{Proof of Lemma 2.1}}

\subsubsection*{Continuity of $\phi$}

By the inequality, for all
$a_1,a_2,a_3,b_1,b_2,b_3\geq0$,
\[
|{\max}\{a_1,a_2,a_3\}-\max\{b_1,b_2,b_3\}|\leq
|a_1-b_1|+|a_2-b_2|+|a_3-b_3|,
\]
we get
%
%
\begin{equation}\label{eq:max0}
|\phi(\alpha_1,\alpha_2,\alpha_3)-\phi(\beta_1,\beta_2,\beta_3)|
\leq\sum_{l=1}^{3} |\alpha_l(c_l)-\beta_l(c_l)|.
\end{equation}
Since $c$ is continuous, if the sequence $((\alpha^n_1,
\alpha^n_2,\alpha^n_3))_{n\geq1}\in\mathcal{M}_b(\mathbb{T})^3$
converges to $(\beta_1,\beta_2,\beta_3)$ (with respect to the
product weak topology), then
\[
{\lim_{n\to\infty}}|\alpha^n_1(c_1)-\beta_1(c_1)|=0,\qquad {\lim
_{n\to\infty}}|\alpha^n_2(
c_2)-\beta_2(c_2)|=0
\]
and
\[
{\lim_{n\to\infty}}|\alpha^n_3(c_3)-\beta_3(c_3)|=0.
\]
The conclusion follows combining these latter three
limits with (\ref{eq:max0}).

\subsubsection*{Continuity of $\Psi$}

For each $l\in\{1,2,3\}$, the projection mapping $\alpha\mapsto
\alpha_{|\mathbb{T}_l}$ is continuous. Hence, the continuity of $\Psi$
follows by the continuity of $\phi$.

\subsubsection*{Continuity of $\Phi$}

Note that, for each fixed $\alpha\in\mathcal{M}_{1}(\mathbb{T})$, it
holds
\[
\Phi(\alpha)=\phi(\alpha_1,\alpha_2,\alpha_3) \qquad\mbox{for some
}\alpha_1,\alpha_2,\alpha_3\in\mathcal{M}_{b}(\mathbb{T})\dvtx
\alpha
_1+\alpha
_2+\alpha_3=\alpha
\]
[indeed, the set
$\{(\alpha_1,\alpha_2,\alpha_3)\in\mathcal{M}_{b}(\mathbb
{T})^3\dvtx\alpha
_1+\alpha_2+\alpha_3=\alpha\}$
is compact with respect to the product weak topology and the
functional $\phi$ is continuous]. For each integer $K>0$, consider
the open covering of $\mathbb{T}$ given by the family formed by the open
balls centered at $x\in\mathbb{T}$ with radius $1/K$. Then by a
classical result (see, e.g., Proposition 16, page 200, in
Royden \cite{royden}) there exists a finite collection
$\{\psi_{n}\}_{1 \leq n \leq N}$ of continuous functions from
$\mathbb{T}$ to $\mathbb{T}$ such that
\begin{eqnarray*}
\sum_{n=1}^{N}\psi_n(x)&=&1 \qquad\mbox{for each }x\in
\mathbb{T},\\
\ell(\operatorname{supp}(\psi_n))&\leq&\pi/K^2 \qquad\mbox{for each
}n=1,\ldots, N.
\end{eqnarray*}
Here the symbol $\operatorname{supp}(\psi_n)$ denotes the support of
$\psi_n$. Let $f$ be a continuous function on $\mathbb{T}$, consider the
modulus of continuity of $f$ defined by
$w_{\delta}(f)={\sup_{|s-t|\leq\delta}}|f(s)-f(t)|$ and set
$f_n=\sup_{x\in \operatorname{supp}(\psi_n)}f(x)$. Note that, for all measures
$\mu\in\mathcal{M}_b(\mathbb{T})$,
%
%
\begin{equation}\label{eq:dom1}\quad
\sum_{n=1}^{N}|\mu(f\psi_n)-f_n\mu(\psi_n)|\leq
w_{{2}/{K}}(f)\sum_{n=1}^{N}\mu(\psi_n)=w_{2/K}(f)
\mu(\mathbb{T}).
\end{equation}
For $i=1,2,3$, define
$r_n^i=\frac{\alpha_i(\psi_n)}{\alpha(\psi_n)}$ if
$\alpha(\psi_n)>0$ and $r_n^i=0$ otherwise. Moreover, for
$\beta\in\mathcal{M}_{b}(\mathbb{T})$, set
%
%
\begin{equation}\label{eq:beta1}
\beta_i(dx)=\sum_{n=1}^{N}r_n^i\psi_n(x)\beta(dx),\qquad i=1,2,3.
\end{equation}
Since
$\alpha_1(\psi_n)+\alpha_2(\psi_n)+\alpha_3(\psi_n)=
\alpha(\psi_n)$, by the properties of the sequence
$\{\psi_n\}_{1\leq n\leq N}$ we have
$\beta_1+\beta_2+\beta_3=\beta$. For any continuous function $f$
on $\mathbb{T}$ we have, for $i=1,2,3$,
%
%
\begin{eqnarray}\label{eq:ineqcont1}
&&
|\beta_i(f)-\alpha_i(f)|\nonumber\\
&&\qquad=\Biggl|\sum_{n=1}^{N} \bigl(r_n^i\beta
(f\psi
_n)-\alpha_i(f\psi_n) \bigr)
\Biggr|\nonumber\\[-8pt]\\[-8pt]
&&\qquad\leq\Biggl|\sum_{n=1}^{N}r_n^i \bigl(\beta(f\psi_n)-\alpha
(f\psi_n) \bigr)\Biggr|+\Biggl|\sum_{n=1}^{N}r_n^i \bigl(f_{n}\alpha
(\psi_n)-\alpha(f\psi_n) \bigr)\Biggr|\nonumber\\
&&\qquad\quad{} +\Biggl|\sum_{n=1}^{N} \bigl(r_n^i
f_n\alpha(\psi_n)-\alpha_i(f\psi_n) \bigr)\Biggr|.\nonumber
\end{eqnarray}
Note that $r_n^i\leq1$, and therefore
%
%
\begin{equation}\label{eq:ineqcont2}
\Biggl|\sum_{n=1}^{N}r_n^i \bigl(\beta(f\psi_n)-\alpha
(f\psi_n) \bigr)\Biggr|\leq {N\max_{1\leq n\leq N}}|\beta
(f\psi_n)-\alpha(f\psi_n)|.
\end{equation}
Using again that $r_n^i\leq1$ and (\ref{eq:dom1}) with
$\mu=\alpha$, we have
%
%
\begin{equation}\label{eq:ineqcont3}\qquad
\Biggl|\sum_{n=1}^{N}r_n^i \bigl(f_n\alpha(\psi_n)-\alpha
(f\psi_n) \bigr)\Biggr|\leq\sum_{n=1}^{N}|f_n\alpha(\psi
_n)-\alpha
(f\psi_n)|\leq w_{2/K}(f).
\end{equation}
By the definition of $r_n^i$ and (\ref{eq:dom1}) it
follows that
%
%
\begin{eqnarray}\label{eq:ineqcont4}
\Biggl|\sum_{n=1}^{N} \bigl(r_n^i
f_n\alpha(\psi_n)-\alpha_i(f\psi_n) \bigr)\Biggr|&=&\Biggl|\sum
_{n=1}^{N} \bigl(
f_n\alpha_i(\psi_n)-\alpha_i(f\psi_n)
\bigr)\Biggr|\nonumber\\[-8pt]\\[-8pt]
&\leq& w_{2/K}(f).\nonumber
\end{eqnarray}
Collecting (\ref{eq:ineqcont1}), (\ref{eq:ineqcont2}),
(\ref{eq:ineqcont3}) and (\ref{eq:ineqcont4}) we have
%
%
\begin{equation}\label{eq:ineqcont5}
|\beta_i(f)-\alpha_i(f)|\leq {N\max_{1\leq n\leq N}}|\beta
(f\psi_n)-\alpha(f\psi_n)|+2w_{2/K}(f).
\end{equation}
Now, let $\{\beta^m\}\subset\mathcal{M}_1(\mathbb{T})$ be a
sequence of probability measures converging to $\alpha$ for the
topology of the weak convergence. We shall prove
\[
\lim_{m\to\infty}\Phi(\beta^m)=\Phi(\alpha).
\]
We first prove
%
%
\begin{equation}\label{eq:upperbound}
\limsup_{m\to\infty}\Phi(\beta^m)\leq\Phi(\alpha).
\end{equation}
Let $K$ be as above and define the Borel measure
$\beta_i^{m}$ as in (\ref{eq:beta1}), with $\beta^m$ in place of
$\beta$ (the definition of $r_n^i$ remains unchanged). By
inequality (\ref{eq:ineqcont5}) and the weak convergence of
$\beta^m$ to $\alpha$, it follows that
\[
{\limsup_{m\to\infty}}|\beta_i ^{m}(f)-\alpha_i(f)|\leq2 w_{2/K} (f).
\]
Applying the above inequality for $f=c_1$, $f=c_2$,
$f=c_3$ and using the inequality (\ref{eq:max0}), we get
\[
{\limsup_{m\to\infty}}|\phi(\beta^m_1,\beta^m_2,\beta^m_3)
-\phi(\alpha_1,\alpha_2,\alpha_3)|\leq6w_{2/K}(c).
\]
Note that by the definition of $\Phi$ and the choice of the
$\alpha_i$'s, $\Phi(\alpha)=\phi(\alpha_1,\alpha_2,\alpha_3)$ and
$\Phi(\beta^m )\leq\phi(\beta^m_1,\beta^m_2,\beta^m_3)$, therefore
\[
\limsup_{m\to\infty}\Phi(\beta^m)\leq\Phi(\alpha)+6w_{2/K}
(c).
\]
The above inequality holds for all $K$, and letting $K$ tend to
infinity, we obtain (\ref{eq:upperbound}). We finally check the
lower semi-continuity bound
%
%
\begin{equation}\label{eq:lowerbound}
\liminf_{m\to\infty}\Phi(\beta^m)\geq\Phi(\alpha).
\end{equation}
Arguing as at the beginning of the proof, we have, for
each fixed $m\geq1$,
\begin{eqnarray}
\Phi(\beta^m)=\phi(\beta_1^m,\beta_2^m,\beta_3^m) \hspace*{50pt}\nonumber\\
\eqntext{\mbox{for
some }
\beta_1^m,\beta_2^m,\beta_3^m\in\mathcal{M}_{b}(\mathbb{T})\dvtx
\beta
_1^m+\beta
_2^m+\beta_3^m=\beta^m.}
\end{eqnarray}
Now, consider an extracted subsequence $(m_k)_{k\geq1}$ such that
\[
\liminf_{m\to\infty}\Phi(\beta^m)=\lim_{k\to\infty}
\phi(\beta_1^{m_k},\beta_2^{m_k},\beta_3^{m_k}).
\]
As already pointed out, $\mathcal{M}_b(\mathbb{T})^3$ is compact with
respect to
the product weak topology. Therefore, up to extracting a
subsequence of $(m_k)_{k\geq1}$, we may assume that
$(\beta_1^{m_k},\beta_2^{m_k},\beta_3^{m_k})$ converges to
$(\beta_1, \beta_2, \beta_3) \in\mathcal{M}_b ( \mathbb{T})^3$.
By construction,
$\beta_1^{m}+\beta_2^{m}+\beta_3^{m}=\beta^m$ and $\beta^m$
converges to $\alpha$, and thus we have
$\beta_{1}+\beta_{2}+\beta_{3}=\alpha$. Then the definition of
$\Phi$ gives
\[
\phi(\beta_{1},\beta_{2},\beta_3)\geq\Phi(\alpha).
\]
Also the continuity of $\phi$ implies
\[
\lim_{k\to\infty}\phi(\beta_1^{m_k},\beta_2^{m_k},\beta_3^{m_k})=
\phi(\beta_1,\beta_2,\beta_3).
\]
The matching lower bound (\ref{eq:lowerbound}) follows.

\subsection{\texorpdfstring{Proof of Lemma \protect\ref{le:tilderho}}{Proof of Lemma 2.2}}

\subsubsection*{Proof of \textup{(i)}}

For each $\alpha\in\mathcal{M}_{1}(\mathbb{T})$, the set
\[
\{(\alpha_1,\alpha_2,\alpha_3)\in\mathcal{M}_{b}(\mathbb
{T})^3\dvtx\alpha
_1+\alpha
_2+\alpha_3
=\alpha\}
\]
is convex; moreover,
the functional $\phi$ is convex on $\mathcal{M}_{b}(\mathbb{T})^3$.
Therefore, by
a classical result of convex analysis, there exists,
$(\alpha_1,\alpha_2,\alpha_3)\in\mathcal{M}_{b}(\mathbb{T})^3$,
such that
$\Phi(\alpha)=\phi(\alpha_1,\alpha_2,\alpha_3)$.

In order to prove that
$\alpha_1(c_1)=\alpha_2(c_2)=\alpha_3(c_3)$, we reason by
contradiction. Assume, for example, that $\Phi(\alpha)=\alpha_1
(c_1)>\max(\alpha_2(c_2),\alpha_3(c_3))$. For $p \in(0,1)$,
define $(\beta_1,\beta_2,\beta_3)=(p\alpha_1,(1-p)\alpha_{1}+
\alpha_2,\alpha_3)$. We have $\beta_{1}+\beta_{2}+\beta_{3}=
\alpha$ and
\[
\phi(\beta_1,\beta_2,\beta_3)=\max\bigl(p\alpha_{1}(c_1),(1-p)\alpha_{1}
(c_2)+\alpha_{2}(c_2),\alpha_{3}(c_3)\bigr).
\]
In particular, for $p$ large enough,
$\phi(\beta_1,\beta_2,\beta_3)=p\alpha_{1}(c_1)<
\phi(\alpha_1,\alpha_2,\alpha_3)$. This is in contradiction with
$\Phi(\alpha)=\phi(\alpha_1,\alpha_2,\alpha_3)$. Now, assume, for
example, that $\Phi(\alpha)=\alpha_{1}(c_1)=\alpha_{2}(c_2)
>\alpha_{3}(c_3)$. The same argument carries over, by considering,
for $p\in(0,1)$,
$(\beta_1,\beta_2,\beta_3)=(p\alpha_1,p\alpha_2,\alpha_{3}+
(1-p)(\alpha_{1}+\alpha_{3}))$. All the remaining cases can be
proved similarly.

\subsubsection*{Proof of \textup{(ii)}}

Since $\mathcal{A}_n\subset\mathcal{B}_n$, we have
$\widetilde\rho_n\leq\rho_n$, and therefore we only need to
establish the claimed lower bound on $\widetilde\rho_{n}$. Let
$B^*$ be an optimal allocation matrix for $\widetilde\rho_n$ and
define the set
\[
I=\bigl\{k\in\{1,\ldots,n\}\dvtx\mbox{there exists }l\in\{1,2,3\}\mbox
{ such
that }b_{kl}^{*}\in(0,1)\bigr\}.
\]
Define the matrix $A=(a_{kl})\in\mathcal{A}_n$ by setting
$a_{kl}=b_{kl}^*$, for any $l\in\{1,2,3\}$, if $k\notin I$, and
$a_{k1}=1$, $a_{k2}=a_{k3}=0$ if $k\in I$. Letting $|I|$ denote
the cardinality of $I$, we have
\begin{eqnarray*}
\widetilde{\rho}_n&=&\max_{1\leq l\leq3} \biggl(\sum_{k\in
I}b_{kl}^*c_l(X_k)+\sum_{k\notin I}b_{kl}^*c_l(X_k) \biggr)\\
&\geq&\max\biggl(\sum_{k\in I}a_{k1}c(X_k)+\sum_{k\notin
I}a_{k1}c(X_k)-|I|\|c\|_{\infty},\max_{l\in\{ 2,3\}} \biggl(\sum_{k
\notin I}a_{kl}c_l(X_k) \biggr) \biggr)\\
&\geq&\max_{1\leq l\leq
3} \biggl(\sum_{k=1}^{n}a_{kl}c_l(X_k) \biggr)-|I|\|c\|_{\infty}
\geq\rho_n-|I|\|c\|_\infty.
\end{eqnarray*}
%
Thus, the claim follows if we prove that $|I|\leq
3$. Reasoning by contradiction, assume that $|I|\geq4$ and, for
$j=1,2,3,4$, denote by $k_j\in I$ four distinct indices in~$I$.
For each $k_j$ there exists $l_j\in\{1,2,3\}$ such that
$b_{k_jl_j}^*\in(0,1)$. Since
\[
b_{k_jl_j}^{*}+\sum_{m\in\{1,2,3\}\setminus\{l_j\}}b_{k_{j}m}^*=1
\]
we deduce that there exist $m_j\in\{1,2,3\}\setminus\{l_j\}$ such
that $b^*_{k_{j}m_j}\in(0,1)$. Thus if $|I|\geq4$, there exist
distinct $k_i,k_j\in\{1,\ldots,n\}$, distinct
$l_i,m_i\in\{1,2,3\}$ and distinct $l_j,m_j\in\{1,2,3\}$ such that
$b_{k_{i}l_i},b_{k_{i}m_i},b_{k_{j}l_j},b_{k_{j}m_j}\in
(0,1)$. Choose $\varepsilon\in
(0,\min\{
b_{k_{i}l_{i}}^*,b_{k_{i}m_{i}}^*,b_{k_{j}l_{j}}^*,b_{k_{j}m_{j}}^*\})$
and define the matrix
$B^{\varepsilon}=(b_{kl}^\varepsilon)\in\mathcal{B}_n$ by
\begin{eqnarray*}
b_{k_{i}l_i}^\varepsilon&=&b_{k_{i}l_i}^{*}-\varepsilon,\qquad
b_{k_{i}m_{i}}^\varepsilon=b_{k_{i}m_i}^{*}+\varepsilon,\\
b_{k_{j}l_j}^\varepsilon&=&b_{k_{j}l_j}^{*}+\varepsilon,\qquad
b_{k_{j}m_j}^\varepsilon=b_{k_{j}m_j}^{*}-\varepsilon,
\end{eqnarray*}
and $b_{kl}^\varepsilon=b_{kl}^*$ otherwise. We define similarly
$B^{-\varepsilon}$ by replacing $\varepsilon$ by
$-\varepsilon$. By part (i) of the lemma, the optimal allocation
matrix $B^{*}$ satisfies
\begin{eqnarray*}
&&\max_{1\leq l,m\leq3} \Biggl(\sum_{k=1}^{n}b_{kl}^{\pm
\varepsilon}c_{l}(X_k),\sum_{k=1}^{n}b_{km}^{\pm
\varepsilon}c_m(X_k) \Biggr)\\
&&\qquad\geq\sum
_{k=1}^{n}b_{k1}^{*}c_1(X_k)
=\sum_{k= 1}^{n}b_{k2}^{*}c_2(X_k)
\\
&&\qquad=\sum_{k=
1}^{n}b_{k3}^{*}c_3(X_k).
\end{eqnarray*}
Therefore
\begin{eqnarray*}
&&\max_{1\leq l,m\leq3} \Biggl(\sum_{k=1}^{n}(b_{kl}^{\pm
\varepsilon}-b_{kl}^{*})c_{l}(X_k),\sum_{k=1}^{n}(b_{km}^{\pm
\varepsilon}-b_{km}^{*})c_m(X_k) \Biggr)\\
&&\qquad =\max\bigl(\mp
\varepsilon\bigl(c_{l_i}(X_{k_i})-c_{l_j}(X_{k_j})\bigr),\pm
\varepsilon\bigl(c_{m_i}(X_{k_i})-c_{m_j}(X_{k_j})\bigr) \bigr)\\
&&\qquad\geq0.
\end{eqnarray*}
It gives $c_{l_i}(X_{k_i})=c_{l_j}(X_{k_j})$ and
$c_{m_i}(X_{k_i})=c_{m_j}(X_{k_j})$ but it a.s. cannot happen
since, by assumption, $\ell(c^{-1}(\{t\}))=0$ for all $t\geq0$.

\subsubsection*{Proof of \textup{(iii)}}

It is an immediate consequence of (ii).

\subsection{A particular cost function: The inverse of signal to noise
plus interference ratio}

\begin{figure}

\includegraphics{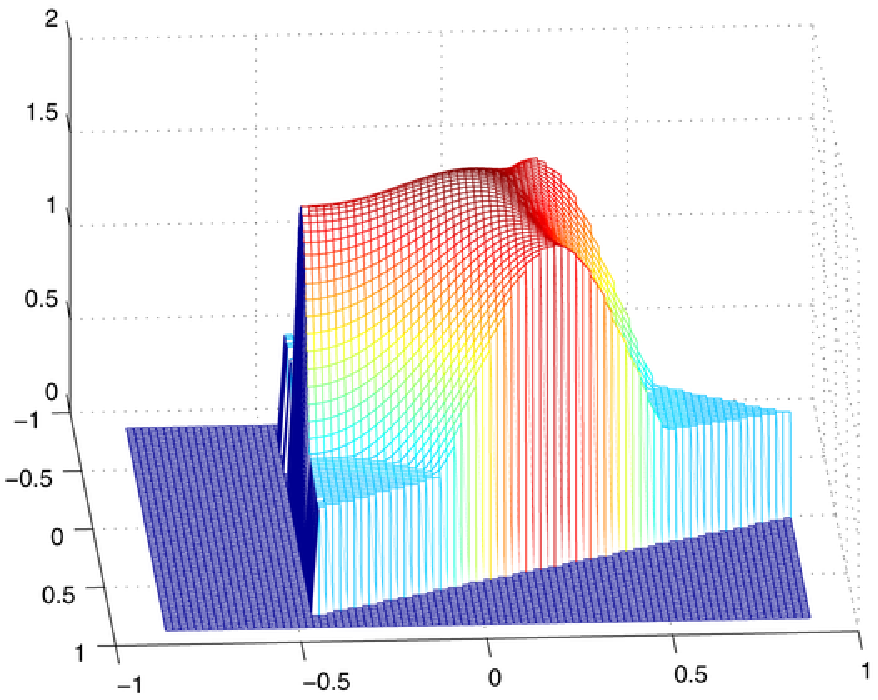}

\caption{The function $L$ with $\alpha=2.5$,
$a=1$ and $b=10$.}\label{fig:cost}
\end{figure}

In this subsection, we prove that the following cost function:
\[
c(x)=
\frac{a+\min\{b,|x-B_2|^{-\alpha}\}+\min\{b,|x-B_3|^{-\alpha}\}
}{\min\{b,|x
-B_1|^{-\alpha}\}},\qquad x\in\mathbb{T},
\]
where $\alpha\geq2$, $a>0$ and $b>(\lambda\sqrt{3}/2)^{-\alpha}$,
satisfies (\ref{ass3}), (\ref{ass4}), (\ref{ass20}), (\ref{ass2})
and (\ref{eq:concave}). To avoid lengthy computations we only checked
numerically the first inequality in (\ref{eq:concave1}). The
typical shape of the function
\[
L(x)=\frac{c_1(x)c_2(x)c_3(x)}{c_1(x)c_2(x)+c_1(x)c_3(
x)+c_2(x)c_3(x)}
\]
is plotted in Figure \ref{fig:cost}, which shows that $L$
attains the supremum at $x=0$. Finally, we show that, for fixed
$\alpha>2$ and $a>0$, for all $b$ large enough, the second inequality in
(\ref{eq:concave1}) holds.

%
%

We first check assumption (\ref{ass3}). We consider only the case
$l=2$, being the case $l=3$ similar. Let $x\in\mathbb{T}$ be such that
$|x-B_1|<|x-B_2|$. Then necessarily, $|x-B_2|>\lambda\sqrt{3}/2$.
With our choice of $b$, we deduce that
\[
\min\{b,|x-B_2|^{-\alpha}\}=|x-B_2|^{-\alpha}<\min\{b,|x
-B_1|^{-\alpha}\}.
\]
By construction
\[
c_2(x)=\frac{a+\min\{b,|x-B_1|^{-\alpha}\}+\min\{b,|x-B_3|^{-\alpha
}\}
}{\min\{b,|x
-B_2|^{-\alpha}\}},\qquad x\in\mathbb{T},
\]
and so (\ref{ass3}) follows easily.

It is immediate to check that $c$ is a Lipschitz function, and the
axial symmetry around the straight line determined by $0$ and
$B_1$ maps $B_2$ into $B_3$. Thus assumptions (\ref{ass4}) and
(\ref{ass2}) follow.

In order to check (\ref{eq:concave}), we note that if $x\in\mathbb{T}_1$,
then, for $l=2,3$, $|x-B_l|\geq|x-B_1|$. Thus, for $l=2,3$,
$\min\{b,|x-B_l|^{-\alpha}\}\leq\min\{b,|x-B_1|^{-\alpha}\}$, and
we deduce
\begin{eqnarray*}
c(x)& = &
\frac{a+\min\{b,|x-B_2|^{-\alpha}\}+\min\{b,|x-B_3|^{-\alpha}\}
}{\min\{b,|x
-B_1|^{-\alpha}\}}\\
& \leq& \frac{a}{ \min\{b,|x
-B_1|^{-\alpha}\}} + 2 \\
& \leq& \lambda^\alpha a + 2 = c(0),
\end{eqnarray*}
where the last inequality is strict if $x\ne0$. Similarly,
$a+\min\{b,|x-B_2|^{-\alpha}\}+\min\{b,|x-B_3|^{-\alpha}\}$ is
minimized for $x=B_1$ and $\min\{b,|x -B_1|^{-\alpha}\}$ is
maximized for $x=B_1$. So, for $x\ne B_1$, $c(x)>c(B_1)$.

Now we check assumption (\ref{ass20}). Define
\[
A_{l}=\{x\in\mathbb{T}\dvtx|x-B_l|<b^{-1/\alpha}\},\qquad l=1,2,3.
\]
With our choice of $b$, if $l\ne m$, we have $A_{l}\cap A_{m}=
\varnothing$. Define
\[
A_0=\mathbb{T}\setminus(A_{1}\cup A_{2}\cup A_{3} ).
\]
Note that, by construction, on each set $A_l$,
$l=0,1,2,3$, the sign of $b-|x-B_m|^{-\alpha}$ is constant for
each $m=1,2,3$. To prove (\ref{ass20}), we shall check that, for
all $t\geq0$ and $l=0,1,2,3$,
%
%
\begin{equation}\label{eq:Al}
\ell\bigl(A_{l}\cap c^{-1}(\{t\}) \bigr)=0.
\end{equation}
We shall only prove the above equality for $l=0$, the other cases
can be shown similarly. Note that
\[
c(x)=|x-B_1|^{\alpha}(a+|x-B_2|^{-\alpha}+|x-B_3|^{-\alpha})\qquad
\forall
x\in A_0.
\]
Using polar coordinates we have
\[
\ell\bigl(A_{0}\cap
c^{-1}(\{t\}) \bigr)=\int_{0}^{2\pi}d\theta\int_{0}
^{\infty}\mathbh{1}\{r\rme^{i\theta}\in
A_0\}\mathbh{1}\{c(r\rme^{i\theta})=t\}r \,dr.
\]
We shall check that, for an arbitrarily fixed $\theta\in
[0,2\pi)$, the function
\[
c_{\theta}(r)=a|r\rme^{i\theta}-B_1|^{\alpha}+
\biggl(\frac{|r\rme^{i\theta}-B_1|}{|r\rme^{i
\theta}-B_2|} \biggr)^{\alpha}+ \biggl(\frac{|r\rme^{i\theta}
-B_1|}{|r\rme^{i\theta}-B_3|} \biggr)^{\alpha},\qquad
r\in I_\theta,
\]
is strictly monotone, where
\[
I_{\theta}=\{r\dvtx r\geq0, r\rme^{i\theta}\in\mathbb{T}\}.
\]
So, for any fixed $\theta\in[0,2\pi)$, the function
$\mathbh{1}\{r\rme^{i\theta}\in A_0\}\mathbh{1}
\{c(r\rme^{i\theta})=t\}$ is different from $0$ for at most
one $r$, and therefore equality (\ref{eq:Al}) for $l =0$
follows. In the following we shall only prove that $c_\theta$ is
strictly decreasing on $I_\theta$ for $\theta\in[-\pi/6,\pi/6]$,
the other cases can be treated similarly. First, note that since
$\theta\in[-\pi/6,\pi/6]$, as $r$ increases, $|r\rme^{i
\theta}-B_1|^{\alpha}$ decreases, while $|r\rme^{i\theta}
-B_3|^{\alpha}$ increases. Thus, $r\mapsto a|r\rme^{i
\theta}-B_1|^{\alpha}$ and
$r\mapsto(\frac{|r\rme^{i\theta}
-B_1|}{|r\rme^{i\theta}-B_3|} )^{\alpha}$ are
decreasing. Note also that, for $\theta\in[-\pi/6,0]$, as $r$
increases, $|r\rme^{i\theta}-B_2|^{\alpha}$ increases. Thus
it suffices to prove that, for a fixed $\theta\in(0,\pi/6]$, the
function
\[
L_\theta(r)=\frac{|r\rme^{i\theta}-B_1|^2}{|r\rme^{i
\theta}-B_2|^2},\qquad r\in\biggl[0,\lambda\biggl(2\cos
\biggl(\frac
{\pi}{6}
-\theta\biggr) \biggr)^{-1} \biggr],
\]
is nonincreasing. Consider the orthonormal basis
$\{\mathbf{e}_1,\mathbf{e}_2\}$ with
$\mathbf{e}_1=\rme^{i{\pi}/{6}}$ and
$\mathbf{e}_{2}=\rme^{-i{\pi}/{3}}$. Setting $\beta=\pi/6
-\theta\in[0,\pi/6)$, $y_{1}=\lambda/2$ and
$y_{2}=\lambda\sqrt{3}/2$, we have
\[
r\rme^{i\theta}=r\cos\beta\mathbf{e}_{1}+r\sin\beta\mathbf{e}_2,\qquad
B_{1}=y_{1}\mathbf{e}_{1}+y_{2}\mathbf{e}_2,\qquad
B_{2}=y_1\mathbf{e}_1-y_2\mathbf{e}_2
\]
and
\[
L_\theta(r)=\frac{(y_{1}-r\cos\beta)^{2}+(y_{2}-r\sin\beta)^{2}
}{(y_{1}-r\cos\beta)^{2}+(y_{2}+r\sin\beta)^2}.
\]
The derivative $L'_\theta(r)$ of $L_\theta(r)$ has the same sign
of
\begin{eqnarray*}
&& - \bigl( \cos\beta( y_1 - r \cos\beta) + \sin\beta( y_2 - r
\sin
\beta) \bigr)
\bigl( (y_1 - r \cos\beta) ^2 + (y_2 + r \sin\beta) ^2 \bigr) \\
&&\qquad{} + \bigl( \cos\beta( y_1 - r \cos\beta) - \sin\beta(
y_2 + r \sin\beta) \bigr)\\
&&\qquad\quad{}\times \bigl( (y_1 - r \cos\beta) ^2 +
(y_2 - r \sin\beta) ^2 \bigr).
\end{eqnarray*}
After simplification, we get easily that $L'_\theta(r)$ has the
same sign of
\[
-2 r \cos\beta\sin\beta- \bigl( ( y_1 - r \cos\beta) ^2 +
y_2 ^2 - r^2 \sin^2 \beta\bigr) \sin\beta.
\]
This last expression is less than or equal to $0$. Indeed, for $r
\in[0,\lambda( 2 \cos\beta)^{-1}]$, we have $0\leq r\sin\beta
\leq y_2$. Hence $L_\theta$ is nonincreasing on its domain.

Finally, we check that, for fixed $\alpha>2$ and $a>0$, it is
possible to determine $b>(\lambda\sqrt{3}/2)^{-\alpha}$ so that
the second inequality in (\ref{eq:concave1}) holds. We deduce
%
%
\begin{eqnarray}\quad
\label{eq:sinr1}
\int_{\mathbb{T}_2}c(x) \,dx&\geq&\int_{\mathbb{T}_2}
\frac{a+\sum_{l=2}^{3}\min\{b,|x-B_l|^{-\alpha}\}}{(\lambda\sqrt
{3}/2)^{-\alpha}}
\,dx\\
\label{eq:sinr2}
&=&\int_{\mathbb{T}_2}
\frac{a+\min\{b,|x-B_2|^{-\alpha}\}+|x-B_3|^{-\alpha}}{(\lambda
\sqrt
{3}/2)^{-\alpha}}
\,dx\\
\label{eq:sinr3}
&\geq&\frac{a/3}{(\lambda\sqrt{3}/2)^{-\alpha}}+ \frac{\pi
b^{1-(2/\alpha)}/6}{(\lambda\sqrt{3}/2)^{-\alpha}}\nonumber\\[-8pt]\\[-8pt]
&&{}+
\bigl(\lambda\sqrt{3}/2\bigr)^{\alpha}\int_{\mathbb{T}_2}|x-B_3|^{-\alpha}
\,dx.\nonumber
\end{eqnarray}
Here (\ref{eq:sinr1}) and (\ref{eq:sinr2}) follow since on
$\mathbb{T}_2$ we have
$|x-B_l|^{-\alpha}<(\lambda\sqrt{3}/2)^{-\alpha}<b$ for $l=1,3$;
(\ref{eq:sinr3}) is consequence of the inequality
$|x-B_2|^{-\alpha}>b$, for any $x\in A_2\cap\mathbb{T}_2$. The claim
follows noticing that, due to our choice of $\alpha$, $c(0)/3$ is
strictly less than the quantity in (\ref{eq:sinr3}), for $b$ large
enough.
\end{appendix}

%

%
\printaddresses

\end{document}